\documentstyle{amsppt}
\magnification=1200
\baselineskip=18pt
\nologo
\TagsOnRight
\document
\define \ts{\thinspace}
\define \di{\partial}
\define \wt{\widetilde}

\define \ZZ{{\Bbb Z}}
\define \sgn{\text{{\rm sgn}}}
\define \de{\text{\rm det}\ts}

\define \ot{\otimes}
\define \wh{\widehat}

\define \gl{{\frak {gl}}}
\define \U{\operatorname {U}}
\define \Z{\operatorname {Z}}
\define \Mat{{\operatorname {Mat}}}

\define \str{\text{{\rm str\ts}}}

\define \C{\Bbb C}
\define \SS{\Bbb S}

\define \Proof{\noindent {\bf Proof. }}

\heading{\bf FACTORIAL SUPERSYMMETRIC SCHUR FUNCTIONS}
\endheading
\heading{\bf AND SUPER CAPELLI IDENTITIES}
\endheading 
\bigskip
\bigskip
\heading{Alexander Molev}\endheading
\bigskip
\bigskip
\bigskip
\bigskip
\noindent
Centre for Mathematics and its Applications,\newline
Australian National University,\newline
Canberra, ACT 0200, Australia\newline
(e-mail: molev\@pell.anu.edu.au)
\bigskip
\bigskip
\bigskip
\bigskip
\noindent
{\bf Abstract}\newline
A factorial analogue of the supersymmetric Schur functions is introduced.
It is shown that factorial versions of the Jacobi--Trudi and
Sergeev--Pragacz formulae hold. The results are applied to construct
a linear basis in the center of the universal enveloping algebra
for the Lie superalgebra $\gl(m|n)$ and to obtain super-analogues
of the higher Capelli identities.

\newpage
\hfill {\sl Dedicated to Professor
A. A. Kirillov on his 60th birthday}  
\bigskip
\noindent
{\bf 0. Introduction}
\bigskip

For a partition $\lambda$ of length $\leq m$
the factorial Schur function $s_{\lambda}(x|a)$
in variables $x=(x_1,\dots,x_m)$
depending on an arbitrary numerical sequence
$a=(a_i)$, $i\in\ZZ$ may be defined
as follows [20, p.\ts 54]. Let
$$
(x|a)^k=(x-a_1)\cdots (x-a_k)
$$
for each $k\geq 0$. Then 
$$
s_{\lambda}(x|a)=
\frac{\de \bigl[(x_j|a)^{\lambda_i+m-i}\bigr]_{1\leq i,j\leq m}}
{\Delta(x)\qquad},\tag 0.1
$$
where $\Delta(x)$ stands for the Vandermonde determinant,
$$
\Delta(x)=\prod_{i<j}(x_i-x_j)=
\de \bigl[(x_j|a)^{m-i}\bigr]_{1\leq i,j\leq m}.
$$
Note that the usual Schur function $s_{\lambda}(x)$ coincides with
$s_{\lambda}(x|a)$ for the zero sequence $a$. The factorial Schur functions
admit many of the classical properties of $s_{\lambda}(x)$ (see
[4--6, 9, 10, 21, 29, 31, 32]), as well as some new ones, e.g.,
the characterization theorem [29] which plays an important role in
the proofs of analogues of the Capelli identity [24, 29].
(The term `factorial' was primarily used only for
the case of the sequence $a$ with $a_i=i-1$;
we use it in the present paper in a broader sense, for an arbitrary
$a$). In particular,
the functions $s_{\lambda}(x|a)$ may be equivalently defined
in terms of tableaux which also enables one to introduce
the skew factorial Schur functions $s_{\lambda/\mu}(x|a)$ for each pair
of partitions $\mu\subset\lambda$.

Super-analogues of the Schur functions can be defined
by using
some specializations of the usual symmetric functions in infinitely many
variables (see [20, p.\ts 58]). This approach goes back to Littlewood [19],
see also [23]. On the other hand, they
naturally emerged in the
representation theory of Lie superalgebras [15, 16] and were studied by
several authors, see, e.g., [2, 3, 7, 9, 13, 33--35, 38--40].

For a pair of partitions $\lambda$ and $\mu$ with
$\mu\subset\lambda$ the supersymmetric skew Schur function
$s_{\lambda/\mu}(x/y)$
in variables $x=(x_1,\dots,x_m)$ and $y=(y_1,\dots,y_n)$
can be defined by the formula
$$
s_{\lambda/\mu}(x/y)=\sum_{\mu\subset\nu\subset\lambda}
s_{\lambda/\nu}(x)s_{\nu'/\mu'}(y),\tag 0.2
$$
where $\lambda'$ denotes the partition conjugate to $\lambda$.

Linear span of the functions (0.2) consists of all supersymmetric
polynomials $P$ in $x$ and $y$; that is, those polynomials which are
symmetric in $x$ and $y$
separately and satisfy the following cancellation property:
the result of setting $x_m=-y_n=z$ in $P$ is independent of $z$.

In this paper we introduce factorial analogues $s_{\lambda/\mu}(x/y\ts|a)$
of the supersymmetric Schur functions parametrized by
arbitrary numerical sequences $a=(a_i)$, $i\in\ZZ$ (see Section 1).
For the zero sequence $a$ this function
coincides with $s_{\lambda/\mu}(x/y)$, and the functions
$s_{\lambda}(x/y\ts|a)$ with $\lambda_{m+1}\leq n$ form a
(non-homogeneous) basis in the
space of supersymmetric polynomials.
We prove that many
properties of the supersymmetric Schur functions have their
factorial analogues.

First, we find the generating series for
the corresponding elementary and complete functions
and check that they are supersymmetric.

Then, using
a modified Gessel--Viennot method [8] we prove an analogue
of the Jacobi-Trudi formula and thus prove that the
polynomials $s_{\lambda/\mu}(x/y\ts|a)$ are also supersymmetric.

Further, we prove a characterization
theorem for the polynomials $s_{\lambda}(x/y\ts|a)$ analogous
to the corresponding theorem for the factorial Schur
polynomials [29] (see also [37]).

Using this theorem we prove
a factorial analogue of the Sergeev--Pragacz formula.
For the usual supersymmetric Schur functions this formula
can be proved by several different ways, see, e.g., [3, 13, 22, 33, 34].
However, these proofs can not be easily carried to the case
of the functions $s_{\lambda}(x/y\ts|a)$,
because for a general sequence $a$ they
lose both the symmetry property of the functions (0.2)
$$
s_{\lambda/\mu}(x/y)=s_{\lambda'/\mu'}(y/x)\tag 0.3
$$
and the specialization property with respect to $x$
$$
s_{\lambda}(x/y)|_{x_m=0}=s_{\lambda}(x'/y),\tag 0.4
$$
where $x'=(x_1,\dots,x_{m-1})$.

In particular, in the case of $a=(0)$ we obtain one more proof of the
Sergeev--Pragacz formula.

As a corollary of the factorial Sergeev--Pragacz formula we get
an analogue of the Berele--Regev factorization theorem for the
polynomials $s_{\lambda}(x/y\ts|a)$ [2], see also [9, 35].

A special case of the factorization theorem yields
an analogue of the dual Cauchy formula which gives a decomposition of
the double product $\prod\prod(x_i+y_j)$ into a sum of products of the
factorial Schur functions. Two other proofs of this formula
are given in [17, 18] and [21].

Goulden and Greene [9] and Macdonald [21] found
a new tableau representation for
the functions
$s_{\lambda/\mu}(x/y)$ in the case of infinite sets of variables
$x=(x_i)$ and $y=(y_i)$, $i\in\ZZ$. 
Using this representation we prove that the corresponding function
$s_{\lambda/\mu}(x/y\ts|a)$ does not depend on $a$ (here $a=(a_i)$ is
regarded as a sequence of independent variables) and coincides with
$s_{\lambda/\mu}(x/y)$.

Finally, we show that many results of the papers [28--31]
concerning higher Capelli identities
and shifted Schur functions
have their natural super-analogues.
In particular, a basis 
in the center of the universal enveloping algebra $\U(\gl(m|n))$
is explicitly constructed. The eigenvalues of the basis elements
in highest weight representations are super-analogues of
the shifted Schur functions.
We outline two proofs of the super-versions of the higher Capelli
identities.
The first proof
uses the characterization theorem for the factorial
supersymmetric Schur polynomials while the second one is based on
the properties of the Jucys--Murphy elements in the group algebra
for the symmetric group.
These identities include
the super Capelli identity found by Nazarov [27].
\medskip

I would like to thank A. Lascoux,
M. Nazarov, A. Okounkov, G. Olshanski\u\i, P. Pragacz for
useful remarks and discussions.

\bigskip
\bigskip
\noindent
{\bf 1. Definition and
combinatorial interpretation of the functions $s_{\lambda/\mu}(x/y\ts|a)$}
\bigskip

We shall suppose that $a=(a_i)$, $i\in\ZZ$ is a fixed sequence of
complex numbers.

For a pair of partitions $\mu\subset\lambda$ the skew factorial
Schur function can be defined by the formula (see, e.g., [9, 21]):
$$
s_{\lambda/\mu}(x|a)=\sum_{T}
\prod_{\alpha\in\lambda/\mu}
(x_{T(\alpha)}-a_{T(\alpha)+c(\alpha)}),\tag 1.1
$$
summed over all semistandard skew tableaux $T$ of shape $\lambda/\mu$ with
entries in the set $\{1,\dots,m\}$, where $T(\alpha)$ is the entry of $T$
in the
cell $\alpha$ and $c(\alpha)=j-i$ 
is the content of $\alpha=(i,j)$. The entries of a semistandard tableau
are supposed to be weakly increasing along rows and strictly increasing
down columns. It can be verified directly (see [9]) that the polynomials
(1.1)
are symmetric in $x$.
It was proved in [21] that for a standard (non-skew) shape $\lambda$
formulae (0.1) and (1.1) define the same function,
which makes the symmetry property obvious for this case.

Due to (1.1), the factorial elementary
and complete symmetric polynomials are given by
$$
e_k(x|a)=s_{(1^k)}(x|a)=\sum_{i_1<\cdots <i_k}(x_{i_1}-a_{i_1})
(x_{i_2}-a_{i_2-1})\cdots (x_{i_k}-a_{i_k-k+1}),\tag 1.2
$$
$$
h_k(x|a)=s_{(k)}(x|a)=
\sum_{i_1\leq\cdots \leq i_k}(x_{i_1}-a_{i_1})
(x_{i_2}-a_{i_2+1})\cdots (x_{i_k}-a_{i_k+k-1}).\tag 1.3
$$

It is immediate from (1.1) that the highest component of
$s_{\lambda/\mu}(x|a)$ is the usual skew Schur polynomial
$s_{\lambda/\mu}(x)$.
This implies that the polynomials $s_{\lambda}(x|a)$ with $l(\lambda)\leq
m$
form a basis in the symmetric polynomials in $x$.
\bigskip
\noindent
{\bf Definition 1.1.}
Let $x=(x_1,\dots,x_m)$ and $y=(y_1,\dots,y_n)$ be two families
of variables. 
Given a sequence $a$
denote by $a^*=(a^*_i)$ another sequence, defined by $a^*_i=-a_{n-i+1}$,
and introduce the
{\it factorial supersymmetric Schur polynomials (functions)} by the
formula:
$$
s_{\lambda/\mu}(x/y\ts|a)=\sum_{\mu\subset\nu\subset\lambda}
s_{\lambda/\nu}(x|a)s_{\nu'/\mu'}(y|a^*).\tag 1.4
$$
\medskip

We shall prove in Section 3 that these polynomials are indeed
supersymmetric,
so their name will be justified.

Comparing (0.2) and (1.4) we see that the highest homogeneous component
of the polynomial
$s_{\lambda/\mu}(x/y\ts|a)$ is $s_{\lambda/\mu}(x/y)$ if
the latter is nonzero, and that
$s_{\lambda/\mu}(x/y\ts|a)$ coincides with $s_{\lambda/\mu}(x/y)$
for the zero sequence $a$.

Using formula (1.1) we can reformulate definition (1.4) in terms of
tableaux.
To distinguish the indices of $x$ and $y$ let us identify
the indices of $y$ with the symbols $1',\dots,n'$.
Consider the diagram of shape $\lambda/\mu$ and fill it with the
indices $1',\dots,n',1,\dots,m$ such that:
\medskip
(a) In each row (resp. column) each primed index is to the left (resp.
above)
from each unprimed index.

(b) Primed indices strictly decrease along rows and weakly decrease
down columns.

(c) Unprimed indices weakly increase along rows and strictly increase
down columns.
\medskip

Denote the resulting tableau by $T$.

\proclaim
{\bf Proposition 1.2} One has the formula
$$
s_{\lambda/\mu}(x/y\ts|a)=\sum_{T}
\prod_{\underset{\ssize T(\alpha)\ts \text{unprimed}}
\to{\alpha\in\lambda/\mu}}
(x_{T(\alpha)}-a_{T(\alpha)+c(\alpha)})
\prod_{\underset{\ssize T(\alpha)\ts \text{primed}}
\to{\alpha\in\lambda/\mu}}
(y_{T(\alpha)}+a_{T(\alpha)+c(\alpha)}).\tag 1.5
$$
\endproclaim

\Proof For each tableau $T$ the cells of the diagram $\lambda/\mu$
occupied by primed indices form a subdiagram $\nu/\mu$ where
$\mu\subset\nu\subset\lambda$. Let us fix such a partition $\nu$ and
sum in (1.5) first over the tableaux $T$
whose primed part forms a subtableau of shape
$\nu/\mu$. The part of such a tableau $T$ formed by unprimed indices
is a semistandard subtableau of shape $\lambda/\nu$ and taking the sum
over these subtableaux we get by (1.1) that
$$
\sum_{T}
\prod_{\alpha\in\lambda/\nu}
(x_{T(\alpha)}-a_{T(\alpha)+c(\alpha)})=s_{\lambda/\nu}(x|a).
$$
It remains to verify that
$$
\sum_{T}\prod_{\alpha\in\nu/\mu}
(y_{T(\alpha)}+a_{T(\alpha)+c(\alpha)})=s_{\nu'/\mu'}(y|a^*),\tag1.6
$$
summed over $\nu/\mu$-tableaux $T$ with entries from $\{1,\dots,n\}$
whose rows strictly decrease and
columns weakly decrease.
Indeed, since $s_{\nu'/\mu'}(y|a^*)$ is symmetric in $y$,
setting $\wt y=(y_n,\dots,y_1)$ and using (1.1) we may write 
$$
s_{\nu'/\mu'}(y|a^*)=s_{\nu'/\mu'}(\wt y\ts|a^*)
=\sum_{T'}\prod_{\alpha'\in\nu'/\mu'}
(\wt y_{T'(\alpha')}-a^*_{T'(\alpha')+c(\alpha')})
$$
$$
=\sum_{T'}\prod_{\alpha'\in\nu'/\mu'}
(y_{n-T'(\alpha')+1}+a_{n-T'(\alpha')+1-c(\alpha')}), \tag 1.7
$$
summed over semistandard $\nu'/\mu'$-tableaux $T'$ with entries from
$\{1,\dots,n\}$. Note that the map
$$
T'(\alpha')\to T(\alpha)=n-T'(\alpha')+1,
$$
where $\alpha=(i,j)\in\nu/\mu$ and $\alpha'=(j,i)\in\nu'/\mu'$, is
a bijection between the set of semistandard $\nu'/\mu'$-tableaux
and the set of $\nu/\mu$-tableaux whose rows strictly decrease and
columns weakly decrease. Obviously, $c(\alpha)=-c(\alpha')$, hence,
(1.7) coincides with the left hand side of (1.6) which completes the
proof.

\bigskip
\bigskip
\noindent
{\bf 2. Generating series for the elementary and complete factorial
supersymmetric polynomials}
\bigskip

Introduce now the {\it elementary\/} and {\it complete factorial
supersymmetric polynomials\/} as special cases of $s_{\lambda}(x/y\ts|a)$
with $\lambda$ being a column or row partition, respectively:
$$
e_k(x/y\ts|a)=s_{(1^k)}(x/y\ts|a)\qquad\text{and}\qquad
h_k(x/y\ts|a)=s_{(k)}(x/y\ts|a).
$$
Using definition (1.4) we can express them in terms of the
polynomials (1.2) and (1.3) as follows:
$$
e_k(x/y\ts|a)=\sum_{p+q=k}e_p(x|\tau^{-q}a)h_q(y|a^*),\tag 2.1
$$
$$
h_k(x/y\ts|a)=\sum_{p+q=k}h_p(x|\tau^qa)e_q(y|a^*),\tag 2.2
$$
where $\tau$ is a shift operator acting on sequences $a$ by
replacing each $a_i$ by $a_{i+1}$. Proposition 1.2 yields the following
explicit formulae for $e_k(x/y\ts|a)$ and $h_k(x/y\ts|a)$:
$$
e_k(x/y\ts|a)=\sum_{p+q=k}
\sum_{\underset{\ssize j_1\geq\cdots\geq j_q}\to
{i_1<\cdots <i_p}}(y_{j_1}+a_{j_1})\cdots (y_{j_q}+a_{j_q-q+1})
(x_{i_1}-a_{i_1-q})\cdots (x_{i_p}-a_{i_p-k+1}),
\tag 2.3
$$
$$
h_k(x/y\ts|a)=\sum_{p+q=k}
\sum_{\underset{\ssize j_1>\cdots > j_q}\to
{i_1\leq\cdots \leq i_p}}(y_{j_1}+a_{j_1})\cdots (y_{j_q}+a_{j_q+q-1})
(x_{i_1}-a_{i_1+q})\cdots (x_{i_p}-a_{i_p+k-1}).
\tag 2.4
$$
We shall suppose that $e_k(x/y\ts|a)=h_k(x/y\ts|a)=0$ if $k<0$.
Note that formulae (2.3) and (2.4) imply the following
symmetry property:
$$
h_k(x/y\ts|a)=e_k(y/x\ts|{-\tau^{k-1}a}).
$$
\bigskip

\proclaim
{\bf Theorem 2.1} One has the following generating series for the
polynomials $e_k(x/y\ts|a)$ and $h_k(x/y\ts|a)$:
$$
1+\sum_{k=1}^{\infty}\frac{(-1)^k e_k(x/y\ts|a)}{(t-a_{m-k+1})\cdots
(t-a_m)}=
\frac{(t-x_1)\cdots (t-x_m)(t-a_1)\cdots (t-a_n)}
{(t-a_1)\cdots (t-a_m)(t+y_1)\cdots (t+y_n)},\tag 2.5
$$
$$
1+\sum_{k=1}^{\infty}\frac{h_k(x/y\ts|a)}{(t-a_{m+1})\cdots (t-a_{m+k})}=
\frac{(t-a_1)\cdots (t-a_m)(t+y_1)\cdots (t+y_n)}
{(t-x_1)\cdots (t-x_m)(t-a_1)\cdots (t-a_n)}.\tag 2.6
$$
\endproclaim

\Proof We shall use induction on $m$. The generating series for the
complete
factorial symmetric polynomials $h_k(y|a)$ is given by
$$
1+\sum_{k=1}^{\infty}\frac{h_k(y|a)}{(t-a_{n+1})\cdots (t-a_{n+k})}=
\frac{(t-a_1)\cdots (t-a_n)}
{(t-y_1)\cdots (t-y_n)}.\tag 2.7
$$
This formula was proved in [31, 32]
in the special case of the sequence $a$ with $a_i=i-1$ and this proof works
in
the general case as well. Note that for $m=0$ we have
$e_k(x/y\ts|a)=h_k(y|a^*)$ by (2.1), and (2.5) follows from (2.7).
Suppose now
that $m\geq 1$. Denote $x'=(x_1,\dots,x_{m-1})$. We see from (2.3) that
$$
e_k(x/y\ts|a)=e_k(x'/y\ts|a)+e_{k-1}(x'/y\ts|a)(x_m-a_{m-k+1}). \tag 2.8
$$
So, using the induction hypotheses, we may write
the right hand side of (2.5) as
$$
\sum_{k=0}^{\infty}\frac{(-1)^k e_k(x'/y\ts|a)}{(t-a_{m-k})\cdots
(t-a_{m-1})}
\cdot\frac{t-x_m}{t-a_m}
$$
$$
=\sum_{k=0}^{\infty}\frac{(-1)^k e_k(x'/y\ts|a)}{(t-a_{m-k+1})\cdots
(t-a_{m})}
\cdot\frac{t-x_m}{t-a_{m-k}}
$$
$$
=\sum_{k=0}^{\infty}\frac{(-1)^k e_k(x'/y\ts|a)+
(-1)^k e_{k-1}(x'/y\ts|a)(x_m-a_{m-k+1})}
{(t-a_{m-k+1})\cdots (t-a_{m})},
$$
where we have used
$$
\frac{t-x_m}{t-a_{m-k}}=1-\frac{x_m-a_{m-k}}{t-a_{m-k}}.
$$
Due to (2.8) this proves (2.5).

The following formula for the generating series for the elementary
factorial symmetric polynomials $e_k(y|a)$ is contained in [20, p.\ts 55]
(see also [31, 32]):
$$
1+\sum_{k=1}^{\infty}\frac{(-1)^k e_k(y|a)}{(t-a_{n-k+1})\cdots (t-a_n)}=
\frac{(t-y_1)\cdots (t-y_n)}
{(t-a_1)\cdots (t-a_n)}.\tag 2.9
$$
For $m=0$ we have $h_k(x/y\ts|a)=e_k(y|a^*)$ by (2.2) and so, for this case
(2.6) follows from (2.9). Now let $m\geq 1$. We see from (2.4) that
$$
h_k(x/y\ts|a)=\sum_{r+s=k}h_r(x'/y\ts|a)(x_m-a_{m+k-s})\cdots
(x_m-a_{m+k-1}).
\tag 2.10
$$
By the induction hypotheses we can write the right hand side of (2.6)
in the form
$$
\sum_{r=0}^{\infty}\frac{h_r(x'/y\ts|a)}{(t-a_{m})\cdots (t-a_{m+r-1})}
\cdot\frac{t-a_m}{t-x_m}
$$
$$
=\sum_{r=0}^{\infty}\frac{h_r(x'/y\ts|a)}{(t-a_{m+1})\cdots (t-a_{m+r})}
\cdot\frac{t-a_{m+r}}{t-x_m}
$$
$$
=\sum_{r=0}^{\infty}\frac{h_r(x'/y\ts|a)}{(t-a_{m+1})\cdots (t-a_{m+r})}
\sum_{s=0}^{\infty}\frac{(x_m-a_{m+r})\cdots (x_m-a_{m+r+s-1})}
{(t-a_{m+r+1})\cdots (t-a_{m+r+s})},
$$
where we have applied (2.7) with $n=1$. The latter expression
can be rewritten as
$$
\sum_{k=0}^{\infty}\frac{1}{(t-a_{m+1})\cdots (t-a_{m+k})}
\sum_{r+s=k}h_r(x'/y\ts|a)(x_m-a_{m+k-s})\cdots (x_m-a_{m+k-1}),
$$
which coincides with the left hand side of (2.6) by (2.10).

\bigskip
\proclaim
{\bf Corollary 2.2} For any $k$
the polynomials $e_k(x/y\ts|a)$ and $h_k(x/y\ts|a)$ are
supersymmetric.
\endproclaim

\Proof Indeed, the cancellation property is obviously satisfied by 
$e_k(x/y\ts|a)$ and $h_k(x/y\ts|a)$ because after setting $x_m=-y_n=z$
the factors on the right hand sides
of (2.5) and (2.6) containing $z$ cancel.
\medskip

\noindent
{\bf Remark.} In his letter to the author A. Lascoux pointed out
that the factorial Schur functions are recovered as
a special case of the double Schubert polynomials [18], [22].
In particular, using the technique of divided differences
one can prove that the complete factorial supersymmetric
polynomials $h_k(x/y\ts|a)$ coincide with
the complete factorial symmetric
polynomials $h_k(x\cup a^{(n)}|{-y}\cup a)$, $a^{(n)}=(a_1,\dots,a_n)$,
and one can also obtain the above generating series for $h_k(x/y\ts|a)$
(see also [1]).

\bigskip
\bigskip
\noindent
{\bf 3. Jacobi--Trudi formula}
\bigskip

The following analogue of the Jacobi--Trudi formula holds for the
functions $s_{\lambda/\mu}(x/y\ts|a)$.

\proclaim
{\bf Theorem 3.1}
$$
s_{\lambda/\mu}(x/y\ts|a)=
\det \bigl[h_{\lambda_i-\mu_j-i+j}
(x/y\ts|\tau^{\mu_j-j+1}a)\bigr]_{1\leq i,j\leq l} \tag 3.1
$$
where $l=l(\lambda)$.
\endproclaim

\Proof We use a modified Gessel--Viennot method [8, 36] (cf.
[9, 31, 32, 35]). Consider a grid consisting of two parts;
the upper half
of the grid if formed by $m$ horizontal lines labelled by $1,\dots,m$
northwards and vertical lines consequently
labelled by the elements of $\ZZ$ eastwards.
For each $i\in\ZZ$  the vertical line labelled by
$i$ breaks into 
two lines at the
intersection point with the horizontal line $1$.
One of the two lines goes south-east and maintain the label $i$ and
the other goes south-west and is labelled by $i'$.

Each vertex of the grid will be denoted by a pair of coordinates $(c,i)$
or $(c,i')$  which are the labels of the lines intersecting at the vertex,
so that $c$ labels a vertical line or a line going south-east.
We shall consider paths
in this grid of the following kind. Each step of a path is north or east
in the upper half of the grid and is north-east or north-west in the
lower part of the grid.
We label each eastern step $(c,i)\to (c+1,i)$ 
of a path with $x_i-a_{i+c}$
and each north-eastern step $(c,i')\to (c+1,i')$ with $y_{i-c}+a_i$.
For a path $\pi$ denote by $L(\pi)$ the product of
these labels with respect to all eastern and north-eastern steps.
Let us check now that
$$
\sum_{\pi} L(\pi)=h_{s-r}(x/y\ts|\tau^{r} a), \tag 3.2
$$
summed over all paths $\pi$ with the initial vertex $(r,(n+r)')$ and the
final vertex $(s,m)$. Indeed, the left hand side of (3.2) can be
calculated as follows. First
fix a number $q$ such that
$0\leq q\leq\text{min}\{s-r,n\}$ and find the sum
$$
\sum_{\pi_{\text{lower}}} L(\pi_{\text{lower}}),\tag 3.3
$$
where $\pi_{\text{lower}}$ runs over
the paths in the lower part of the grid 
with the initial vertex $(r,(n+r)')$ and the
final vertex $(r+q, (r+q)')$
(which belongs to the horizontal line 1 and has also
the coordinates $(r+q,1)$ if regarded as a vertex of the upper half
of the grid).
Suppose that such a path $\pi_{\text{lower}}$
has north-eastern steps of the form
$$
\align
(r,(j_1+r)')&\to (r+1,(j_1+r)'),\\
&\cdots\\
(r+q-1,(j_q+r)')&\to (r+q,(j_q+r)'),
\endalign
$$
where $n\geq j_1\geq\cdots\geq j_q\geq q$. Then the product of
the labels of $\pi_{\text{lower}}$ 
with respect to the north-eastern steps equals
$$
(y_{j_1}+a_{j_1+r})\cdots (y_{j_q-q+1}+a_{j_q+r}).
$$
Hence, the sum (3.3) equals
$$
\align
&\sum_{n\geq j_1\geq\cdots\geq j_q\geq q}
(y_{j_1}+a_{j_1+r})\cdots (y_{j_q-q+1}+a_{j_q+r})\\
=&\sum_{n\geq i_1>\cdots> i_q\geq 1}
(y_{i_1}+a_{i_1+r})\cdots (y_{i_q}+a_{i_q+q-1+r}).
\endalign
$$
Applying (1.6) for $\nu/\mu=(q)$ we conclude that this expression
coincides with $e_{q}(y|(\tau^r a)^*)$.

Similarly, one can easily check that the sum 
$$
\sum_{\pi_{\text{upper}}} L(\pi_{\text{upper}}),
$$
where $\pi_{\text{upper}}$ runs over all paths in the upper half of
the grid with the initial vertex $(r+q,1)$ and the
final vertex $(s,m)$, coincides with $h_{s-r-q}(x|\tau^{r+q}a)$.

Thus, the left hand side of (3.2) equals
$$
\sum_q h_{s-r-q}(x|\tau^{r+q}a)e_{q}(y|(\tau^r a)^*),
$$
which coincides with $h_{s-r}(x/y\ts|\tau^{r} a)$ by (2.2).

A straightforward application of the Gessel--Viennot argument
with the use of (3.2) shows that the determinant on the right hand side of
(3.1) can be represented as
$$
\sum_{\Cal P}L(\pi_1)\cdots L(\pi_l),\tag 3.4
$$
summed over sets of nonintersecting paths
$\Cal P=(\pi_1,\dots,\pi_l)$, where the initial vertex of $\pi_i$
is $(\mu_i-i+1,(n+\mu_i-i+1)')$ and the final vertex is
$(\lambda_i-i+1, m)$.

It remains to prove that (3.4) coincides with $s_{\lambda/\mu}(x/y\ts|a)$.
Given set $\Cal P=(\pi_1,\dots,\pi_l)$ of nonintersecting paths we
construct a $\lambda/\mu$-tableau $T$ as follows. Note that the total
number
of north-eastern and eastern steps of the path $\pi_i$ equals
$\lambda_i-\mu_i$. Let us numerate these steps starting with
the initial vertex of $\pi_i$.
If the $k$th step
is north-eastern of the form $(c,j')\to (c+1,j')$
then the entry of the $k$th cell in the $i$-th row of $\lambda/\mu$
is $(j-c)'$ and if the $k$th step
is eastern of the form $(c,j)\to (c+1,j)$ then the entry of
this cell is $j$. It can be easily seen that this correspondence is
a bijection between the sets of nonintersecting paths and the tableaux
used in the combinatorial interpretation (1.5) of the polynomials
$s_{\lambda/\mu}(x/y\ts|a)$. Moreover, the corresponding summands in
(3.4) and (1.5) are clearly the same,
which proves the theorem. 
\medskip

Note that for $n=0$ formula (3.1) turns into
the Jacobi--Trudi formula for the factorial Schur functions
$s_{\lambda/\mu}(x|a)$ while for $m=0$ we get a factorial analogue
of the N\"agelsbach--Kostka formula (see [20, p.\ts 56]) and setting
$a=(0)$
we obtain the Jacobi--Trudi formula for the functions
$s_{\lambda/\mu}(x/y)$ [35].
\medskip

Corollary 2.2 and Theorem 3.1 imply
\bigskip

\proclaim
{\bf Corollary 3.2} The polynomials $s_{\lambda/\mu}(x/y\ts|a)$
are supersymmetric.
\endproclaim
\bigskip

An analogue of (3.1) for the factorial elementary supersymmetric
polynomials can be proved in the same way by using
a grid obtained from the one we have used in the above proof
by interchanging its upper and lower parts.
We state this formula here
without proof.

\proclaim
{\bf Theorem 3.3}
$$
s_{\lambda'/\mu'}(x/y\ts|a)=
\det \bigl[e_{\lambda_i-\mu_j-i+j}
(x/y\ts|\tau^{-\mu_j+j-1}a)\bigr]_{1\leq i,j\leq l}. \tag 3.5
$$
\endproclaim

\bigskip
\bigskip
\noindent
{\bf 4. Characterization theorem}
\bigskip

Our aim in this section is to obtain characterization properties
for the supersymmetric polynomials (Theorems 4.5 and 4.5$'$)
which will play an important role
in our proof of a factorial analogue of the Sergeev--Pragacz formula
(Theorem 5.1) and super-analogues of the Capelli identity (Theorem 8.1).
We start with investigating vanishing and specialization
properties of the polynomials $s_{\lambda}(x/y\ts|a)$.
\medskip

We say that a partition $\lambda$ is contained in the $(m,n)$-hook if
$\lambda_{m+1}\leq n$. Note that
$$
s_{\lambda}(x/y\ts|a)=0 \quad\text{unless}\quad
\lambda\subset(m,n)\text{-hook}.
\tag 4.1
$$
Indeed, by definition (1.4)
$$
s_{\lambda}(x/y\ts|a)=\sum_{\rho\subset\lambda}
s_{\rho'}(y|a^*)s_{\lambda/\rho}(x|a). \tag 4.2
$$
Clearly, $s_{\rho'}(y|a^*)=0$ unless $\rho_1\leq n$. If $\lambda$ is
not contained in the $(m,n)$-hook then $\lambda'_{n+1}>m$. So, if
$\rho_1\leq n$ then $s_{\lambda/\rho}(x|a)=0$.
\medskip

For each $\lambda\subset(m,n)$-hook we
introduce two partitions $\mu=\mu(\lambda)$ and $\nu=\nu(\lambda)$
as follows.
Nonzero parts of $\mu$ are defined by $\mu_i=\lambda_i-n$ for $\lambda_i>n$
and nonzero parts of $\nu$ are defined by $\nu_j=\lambda'_j-m$
for $\lambda'_j>m$.

For any partition $\alpha$ such that $l(\alpha)\leq l$ introduce
the $l$-tuple
$$
a_{\alpha}=a_{\alpha}^{(l)}=
(a_{\alpha_1+l},\dots,a_{\alpha_l+1})
$$
of elements of the sequence $a$.
We shall only consider $m$- or $n$-tuples $a_{\alpha}$ as values for
the variables $x=(x_1,\dots,x_m)$ or $y=(y_1,\dots,y_n)$ so we
shall not specify the number $l$ if it is clear from the context.

We shall often use the following vanishing properties of the factorial
Schur polynomials $s_{\lambda}(x|a)$ (see [29, 31]). For any
partition $\sigma$ such that $l(\sigma)\leq m$ and
$\lambda\not\subset\sigma$ one has
$$
s_{\lambda}(a_{\sigma}|a)=0.\tag 4.3
$$
Moreover, if $l(\lambda)\leq m$ then
$$
s_{\lambda}(a_{\lambda}|a)=
\prod_{(i,j)\in\lambda}(a_{\lambda_i+m-i+1}-a_{m-\lambda'_j+j}).\tag 4.4
$$
In particular,
$$
s_{\lambda}(a_{\lambda}|a)\ne 0,
$$
provided that the sequence $a$ is multiplicity free, that is,
$a_i\ne a_j$ if $i\ne j$. Relations (4.3) and (4.4) can be
easily deduced from either (0.1), or (1.1). In the special case of the
sequence $a$ with $a_i=i+\text{const}$
the right hand side of (4.4) turns into
the product $H(\lambda)$ of the hook lengths of all cells of $\lambda$
(see [29]).
\bigskip

\proclaim
{\bf Proposition 4.1} Let $\eta$ be a partition of length $\leq n$
such that $\nu\not\subset\eta$. Then
$$
s_{\lambda}(x/a^*_{\eta}\ts|a)=0. \tag 4.5
$$
\endproclaim

\Proof We use (4.2). 
Note that $s_{\lambda/\rho}(x|a)=0$ unless $\lambda'_i-\rho'_i\leq m$
for any $i$. This means that the sum in (4.2) can be only taken
over $\rho\subset\lambda$ such that $\nu'\subset \rho$. By the assumption,
$\nu\not\subset\eta$, hence, $\rho'\not\subset\eta$. But in this case
$s_{\rho'}(a^*_{\eta}\ts|a^*)=0$ by (4.3). Thus, all summands in (4.2)
vanish for $y=a^*_{\eta}$ which completes the proof.
\bigskip

\proclaim
{\bf Proposition 4.2} Let $\gamma$ be a partition of length $\leq m$
such that $\mu\not\subset\gamma$. Then
$$
s_{\lambda}((\tau^n a)_{\gamma}/y\ts|a)=0. \tag 4.6
$$
\endproclaim

\Proof We use again that $s_{\rho'}(y|a^*)=0$ in (4.2) unless $\rho_1\leq
n$.
Relation (4.6) will follow from the fact that
for any such $\rho$ one has 
$$
s_{\lambda/\rho}((\tau^n a)_{\gamma}|a)=0.\tag 4.7
$$
This relation is a simple generalization of the vanishing
property (4.3)
and can be proved by similar arguments; see [29, 31]. Indeed,
since the left hand side of (4.7) is a polynomial in $a$,
we may assume without loss of generality that the sequence
$a$ is multiplicity free.
The polynomials $s_{\lambda/\rho}(x|a)$ are symmetric in
$x$, so replacing $x$ with $\wt x=(x_m,\dots,x_1)$,
we may rewrite definition (1.1) in the following form: 
$$
s_{\lambda/\rho}(x|a)=\sum_{T}
\prod_{\alpha\in\lambda/\rho}
(x_{m-T(\alpha)+1}-a_{T(\alpha)+c(\alpha)}), \tag 4.8
$$
summed over semistandard $\lambda/\rho$-tableau $T$ with entries from
$\{1,\dots,m\}$. We shall verify that all summands in (4.8) vanish
for $x=(\tau^n a)_{\gamma}$. Indeed, let us suppose that for some tableau
$T$
$$
\prod_{\alpha\in\lambda/\rho}
(\{(\tau^n a)_{\gamma}\}_{m-T(\alpha)+1}-a_{T(\alpha)+c(\alpha)})\ne 0
$$
or equivalently,
$$
\prod_{\alpha\in\lambda/\rho}
(a_{\gamma_{m-T(\alpha)+1}+T(\alpha)+n}-a_{T(\alpha)+c(\alpha)})\ne 0.
$$
Since $a$ is multiplicity free, this implies that
$$
\gamma_{m-T(\alpha)+1}+n\ne c(\alpha)\tag 4.9
$$
for all $\alpha$. For the entries of the first row of the tableau $T$
we have
$$
T(1,n+1)\leq\cdots\leq T(1,n+\mu_1). \tag 4.10
$$
Applying (4.9) for $\alpha=(1,n+1)$ we obtain 
$
\gamma_{m-T(1,n+1)+1}\geq 1.
$
Further, by (4.10),
$$
\gamma_{m-T(1,n+2)+1}\geq\gamma_{m-T(1,n+1)+1}\geq 1.
$$
By (4.9), applied for $\alpha=(1,n+2)$ we then have
$
\gamma_{m-T(1,n+2)+1}\geq 2.
$
Similarly, using an easy induction argument we see
that for any $i=1,\dots,\mu_1$
$$
\gamma_{m-T(1,n+i)+1}\geq i.
$$
On the other hand, for the entries of the $(n+i)$th column of $T$ we have
$$
T(1,n+i)<\cdots <T(\mu'_i,n+i).
$$
Hence,
$$
\gamma_{m-T(\mu'_i,n+i)+1}\geq\cdots\geq
\gamma_{m-T(1,n+i)+1}\geq i.
$$
This means that $\gamma'_i\geq\mu'_i$, and so, $\mu\subset\gamma$ which
contradicts to the assumption of the proposition. The proof is complete.
\medskip

As we pointed out in Introduction, for a general sequence $a$
the polynomials $s_{\lambda}(x/y\ts|a)$
lose the specialization property of $s_{\lambda}(x/y)$
with respect to $x$. However, an analogue
of this property with respect to $y$ still holds.
\bigskip

\proclaim
{\bf Proposition 4.3} One has
$$
s_{\lambda}(x/y\ts|a)|_{y_n=-a_n}=s_{\lambda}(x/y'\ts|a), \tag 4.11
$$
where $y'=(y_1,\dots,y_{n-1})$.
\endproclaim

\Proof This follows from the specialization property of the
factorial Schur polynomials (see [31, 32]). Indeed, by (1.6) we have
$$
s_{\rho'}(y|a^*)=\sum_{T}\prod_{\alpha\in\rho}
(y_{T(\alpha)}+a_{T(\alpha)+c(\alpha)}), \tag 4.12
$$
summed over $\rho$-tableaux $T$ with entries from $\{1,\dots,n\}$
whose rows strictly decrease and
columns weakly decrease. If for a tableau $T$ one has $T(1,1)=n$, then
since $c(1,1)=0$ the corresponding summand
in (4.12) vanishes for $y_n=-a_n$. So we have the property:
$$
s_{\rho'}(y|a^*)|_{y_n=-a_n}=s_{\rho'}(y'|{a^*}'),
$$
where ${a^*}'$ is the sequence $({a_i^*}')$ with ${a_i^*}'=-a_{n-i}$.
Now (4.11) follows from (4.2).
\medskip

We can now prove a vanishing theorem for the polynomials
$s_{\lambda}(x/y\ts|a)$ (cf. [29, 31]). Let $\zeta$ be a partition which is
contained in the $(m,n)$-hook. Introduce two other partitions
$\xi=\xi(\zeta)$
and $\eta=\eta(\zeta)$
as follows:
$
\xi=(\zeta_1,\dots,\zeta_m),
$
and nonzero parts of $\eta=(\eta_1,\dots,\eta_n)$ are defined by
$\eta_i=\zeta'_i-m$,
if $\zeta'_i> m$.

In particular, $l(\xi)\leq m$ and $l(\eta)\leq n$ and we may consider
the $m$-tuple $a_{\xi}$ of elements of the sequence $a$ and
the $n$-tuple $a^*_{\eta}$ of elements of the sequence $a^*$.
\bigskip

\proclaim
{\bf Theorem 4.4} Let $\lambda,\zeta$ be partitions which are contained in
the $(m,n)$-hook. 

{\rm (i)} If $\lambda\not\subset\zeta$ then
$$
s_{\lambda}(a_{\xi}/a^*_{\eta}\ts|a)=0. \tag 4.13
$$

{\rm (ii)} If $\lambda=\zeta$ then
$$
s_{\lambda}(a_{\xi}/a^*_{\eta}\ts|a)=
\prod_{(i,j)\in\lambda}(a_{\lambda_i+m-i+1}-a_{m-\lambda'_j+j}).\tag 4.14
$$
In particular, if the sequence $a$ is multiplicity free then
$
s_{\lambda}(a_{\xi}/a^*_{\eta}\ts|a)\ne 0.
$
\endproclaim

\Proof Denote by $r$ the length of the partition $\eta$. Then
$$
a^*_{\eta}=(-a_{1-\eta_1},\dots,-a_{r-\eta_r},-a_{r+1},\dots,-a_n).
$$
By Proposition 4.3, the result of setting $y_i=-a_i$ for $i=r+1,\dots,n$
in the polynomial $s_{\lambda}(x/y\ts|a)$ is the polynomial
$s_{\lambda}(x/y^{(r)}\ts|a)$, where $y^{(r)}=(y_1,\dots,y_r)$.
So, (4.13) and (4.14) can be regarded as relations for the
families of variables $x$ and $y^{(r)}$.
Moreover, due to (4.1) we can only consider the case when
$\lambda$ is contained
in the $(r,m)$-hook. In other words,
we may assume without loss of generality that the length of $\eta$ is $n$.
In particular, the partition $(n^m)$ is contained in $\zeta$.

Now let us prove (i). By Proposition 4.1 we may assume that
$\nu\subset\eta$.
Since $(n^m)\subset\xi$ we may write $a_{\xi}=(\tau^n a)_{\gamma}$, where
the partition $\gamma$ is defined by $\gamma_i=\xi_i-n$. Hence, if
$s_{\lambda}(a_{\xi}/a^*_{\eta}\ts|a)$ is nonzero then
$\mu\subset\gamma$ by Proposition 4.2.
Since $\lambda\subset(m,n)$-hook this implies that
$\lambda\subset\zeta$ and
(i) is proved.

To prove (ii) we note that
since $(n^m)\subset\zeta=\lambda$ we have
$$
a_{\xi}=(\tau^n a)_{\mu}\qquad\text{and}\qquad
a^*_{\eta}=a^*_{\nu}. \tag 4.15
$$
As we have noticed in the proof of Proposition 4.1, the sum
in (4.2) can be only taken over those partitions $\rho\subset\lambda$
for which $\nu'\subset\rho$. On the other hand,
using (4.3)
we find that
$s_{\rho'}(a^*_{\nu}\ts|a^*)=0$ unless $\rho'\subset\nu$.
Hence, for $y=a^*_{\nu}$ relation (4.2) turns into
$$
s_{\lambda}(x/a^*_{\nu}\ts|a)=s_{\nu}(a^*_{\nu}\ts|a^*)\ts
s_{\lambda/\nu'}(x|a).\tag 4.16
$$

Now let us write the polynomial $s_{\lambda/\nu'}(x|a)$
in terms of tableaux using
(1.1).
By definition of $\nu$ we
have $\lambda'_i-\nu_i=m$ for all $i=1,\dots n$. Since the
tableaux $T$ in (1.1)
are column strict, all of them have the same entries in the first
$n$ columns, namely, the numbers $1,\dots,m$ written in each of these
columns downwards. On the other hand, the entries of the
subdiagram $\mu$ can form an arbitrary semistandard $\mu$-tableau.
For a cell $(i,n+j)\in\lambda$ we obviously have
$c(i,n+j)=n+c(i,j)$, so definition (1.1) for $s_{\lambda/\nu'}(x|a)$
takes now the form
$$
s_{\lambda/\nu'}(x|a)=s_{\mu}(x|\tau^n a)
\prod_{i=1}^m\prod_{j=1}^n(x_i-a_{j-\nu_j}).
$$
We then obtain from (4.16) that
$$
s_{\lambda}((\tau^na)_{\mu}/a^*_{\nu}\ts|a)=
s_{\nu}(a^*_{\nu}\ts|a^*)\ts
s_{\mu}((\tau^na)_{\mu}|\tau^n a)
\prod_{i=1}^m\prod_{j=1}^n(a_{\mu_i+m+n-i+1}-a_{j-\nu_j}).\tag 4.17
$$
Now (4.14) follows from (4.4).
The theorem is proved.
\medskip

Theorem 4.4 implies that supersymmetric polynomials in $x$ and $y$
of degree $\leq k$
are characterized by their values at $x=a_{\xi}$ and $y=a^*_{\eta}$,
where $\xi=\xi(\zeta)$ and $\eta=\eta(\zeta)$ with $|\zeta|\leq k$.
More exactly, we have the following result (cf. [29, 31, 37]).
\bigskip

\proclaim
{\bf Theorem 4.5} Suppose that the sequence $a$ is multiplicity free. 
Let $f(x/y)$ and $g(x/y)$ be supersymmetric polynomials of degree $\leq k$
such that
$$
f(a_{\xi}/a^*_{\eta})=g(a_{\xi}/a^*_{\eta}) \tag 4.18
$$
for any partition $\zeta\subset(m,n)$-hook
with $|\zeta|\leq k$. Then $f(x/y)=g(x/y)$.
\endproclaim

\Proof It is well-known that the functions $s_{\lambda}(x/y)$ with
$\lambda\subset(m,n)$-hook form a basis in supersymmetric polynomials
in $x$ and $y$ (see, e.g., [20, p.\ts 61]). So do the functions
$s_{\lambda}(x/y\ts|a)$, because the highest term of
$s_{\lambda}(x/y\ts|a)$
is $s_{\lambda}(x/y)$. Hence, we can write the polynomial
$f(x/y)-g(x/y)$ as a linear combination of the $s_{\lambda}(x/y\ts|a)$:
$$
f(x/y)-g(x/y)=\sum_{\lambda}c_{\lambda}s_{\lambda}(x/y\ts|a).
$$
Moreover, since the degree of the polynomial on the
left hand side $\leq k$ we may assume that $|\lambda|\leq k$.
Introduce any total order on the set of partitions such that
$|\lambda|<|\mu|$ implies $\lambda<\mu$. The condition (4.18)
yields the following homogeneous system of linear equations on the
coefficients $c_{\lambda}$:
$$
\sum_{\lambda}c_{\lambda}s_{\lambda}(a_{\xi}/a^*_{\eta}\ts|a)=0,\qquad
|\lambda|,|\zeta|\leq k.
$$
Theorem 4.4 implies that the matrix
$(s_{\lambda}(a_{\xi}/a^*_{\eta}\ts|a))_{\lambda,\zeta}$
of this system, whose rows
and columns are arranged in accordance with
this order, is triangular with nonzero diagonal
elements. Hence, $c_{\lambda}\equiv 0$ which proves the theorem.
\medskip

Theorem 4.5 can be obviously reformulated in the following equivalent form.
\bigskip

\proclaim
{\bf Theorem 4.5$'$} Suppose that the sequence $a$ is multiplicity free.
Let $f(x/y)$ be a supersymmetric polynomial such that
$$
f(x/y)=s_{\lambda}(x/y)+\text{lower terms}
$$
for some partition $\lambda\subset(m,n)$-hook and
$
f(a_{\xi}/a^*_{\eta})=0
$
for any partition $\zeta\subset(m,n)$-hook with $|\zeta|<|\lambda|$. Then
$f(x/y)=s_{\lambda}(x/y\ts|a)$.
\endproclaim

\bigskip
\bigskip
\noindent
{\bf 5. Factorial Sergeev--Pragacz formula}
\bigskip

In this section we apply Theorem 4.5 for the proof
of an analogue of the Sergeev--Pragacz formula for the
polynomials $s_{\lambda}(x/y\ts|a)$. 
In particular, for $a=(0)$ we get one more
proof of the original formula (cf. [3, 13, 22, 33, 34]).

Suppose a partition $\lambda$ is contained in the $(m,n)$-hook.
Define the partitions $\mu$ and $\nu$ as in Section 4 and denote by
$\rho=(\rho_1,\dots,\rho_m)$
the part of $\lambda$ which is contained in the rectangle
$(n^m)$, that is,
$\rho_i=\text{min}\{\lambda_i,n\}$.

The following
analogue of the Sergeev--Pragacz formula holds.

\proclaim
{\bf Theorem 5.1}
$$
s_{\lambda}(x/y\ts|a)=\frac{\dsize
\sum_{\sigma\in S_m\times S_n}\sgn(\sigma)\cdot
\sigma\bigl\{f_{\lambda}(x/y\ts|a)\bigr\}}{\Delta(x)\Delta(y)},\tag 5.1
$$
where
$$
f_{\lambda}(x/y\ts|a)=(x_1|\tau^{\rho_1}a)^{\mu_1+m-1}\cdots
(x_m|\tau^{\rho_m}a)^{\mu_m}(y_1|a^*)^{\nu_1+n-1}\cdots(y_n|a^*)^{\nu_n}
\prod_{(i,j)\in\rho}(x_i+y_j).
$$
\endproclaim

\Proof First of all we note that both sides of (5.1) depend
polynomially on $a$ so we may assume without loss of
generality that the sequence $a$ is multiplicity free.

Denote the right hand side of (5.1) by $\varphi_{\lambda}(x/y\ts|a)$.
To apply Theorem 4.5 we have to verify that this polynomial
is supersymmetric and that $s_{\lambda}(x/y\ts|a)$ and 
$\varphi_{\lambda}(x/y\ts|a)$
have the same values at $x=a_{\xi}$ and $y=a^*_{\eta}$ for any
$\zeta\subset(m,n)$-hook such that $|\zeta|\leq |\lambda|$.

The polynomials $\varphi_{\lambda}(x/y\ts|a)$ are obviously symmetric
in $x$ and $y$, and so, to prove that they are supersymmetric
we only need to check that they satisfy
the cancellation property. This can be done exactly in the same way as
in the case $a=(0)$ (see e.g., [20, p.\ts 61], [33, 34]).

Let us check now that Propositions 4.1--4.3 hold for
the polynomials $\varphi_{\lambda}(x/y\ts|a)$ too. To check (4.5) we
represent the numerator of the right hand side of (5.1) in the following
form:
$$
\sum_{\sigma\in S_m}\sgn(\sigma)\cdot\sigma
\bigl\{(x_1|\tau^{\rho_1}a)^{\mu_1+m-1}\cdots
(x_m|\tau^{\rho_m}a)^{\mu_m}\ts g_{\lambda}(x/y\ts|a)\bigr\},
$$
where
$$
g_{\lambda}(x/y\ts|a)=
\det\bigl
[(y_j|a^*)^{\nu_i+n-i}(y_j+x_1)\cdots(y_j+x_{\rho'_i})\bigr]_{1\leq i,j\leq
n}.
$$
The condition $\nu\not\subset\eta$ means that there exists $k$ such that
$\eta_k<\nu_k$. For $y=a^*_{\eta}$ the factor
$
(y_j|a^*)^{\nu_i+n-i}
$
takes the value
$$
((a^*_{\eta})_j|a^*)^{\nu_i+n-i}=(a^*_{\eta_j+n-j+1}-a^*_1)\cdots
(a^*_{\eta_j+n-j+1}-a^*_{\nu_i+n-i}).\tag 5.2
$$
On the other hand, if $i\leq k\leq j$ then
$$
1\leq \eta_j+n-j+1\leq \eta_k+n-k+1\leq \nu_k+n-k\leq \nu_i+n-i.
$$
This implies that (5.2) is zero and all the $ij$th entries of the
determinant $g_{\lambda}(x/a^*_{\eta}\ts|a)$ with $i\leq k\leq j$ are zero
and so $g_{\lambda}(x/a^*_{\eta}\ts|a)=0$. Since the Vandermonde
determinant
$\Delta(y)$ does not vanish for $y=a^*_{\eta}$ this proves the assertion.
(This argument is very similar to that used in [29, 31] for the proof of
(4.3)).

To check that the polynomials $\varphi_{\lambda}(x/y\ts|a)$ satisfy (4.6)
we rewrite the numerator of the right hand side of (5.1) in the
form:
$$
\sum_{\sigma\in S_n}\sgn(\sigma)\cdot\sigma
\bigl\{(y_1|a^*)^{\nu_1+n-1}\cdots(y_n|a^*)^{\nu_n}
\ts h_{\lambda}(x/y\ts|a)\bigr\},
$$
where
$$
h_{\lambda}(x/y\ts|a)=
\det\bigl
[(x_j|\tau^{\rho_i}a)^{\mu_i+m-i}(x_j+y_1)\cdots(x_j+y_{\rho_i})
\bigr]_{1\leq i,j\leq m}.
$$
The condition $\mu\not\subset\gamma$ implies that for some $k$ one has
$\gamma_k<\mu_k$. In particular, this means that $\rho_k=n$ and hence
$\rho_1=\cdots=\rho_k=n$. Therefore, the $ij$th entry of the
determinant $h_{\lambda}(x/y\ts|a)$ for $i\leq k\leq j$ has the form
$$
(x_j|\tau^na)^{\mu_i+m-i}(x_j+y_1)\cdots(x_j+y_{\rho_i}).
$$
Repeating the previous argument, we conclude that
$h_{\lambda}((\tau^n a)_{\gamma}/y\ts|a)=0$, which completes the proof.

Let us prove now that
$$
\varphi_{\lambda}(x/y\ts|a)|_{y_n=-a_n}=
\varphi_{\lambda}(x/y'\ts|a), \tag 5.3
$$
where $y'=(y_1,\dots,y_{n-1})$ and we define
$\varphi_{\lambda}(x/y\ts|a)=0$ if $\lambda\not\subset(m,n)$-hook.
Indeed, since $a^*_1=-a_n$, we obviously have for $\nu_n>0$ that
$$
\varphi_{\lambda}(x/y\ts|a)|_{y_n=-a_n}=0,
$$
and hence (5.3) is true, because $\lambda$ is not contained in the
$(m,n-1)$-hook. So, we can suppose that $\nu_n=0$. Since $\nu_i+n-i>0$
for $i=1,\dots,n-1$, after setting $y_n=-a_n$
we may restrict the sum in the numerator
of the right hand side of (5.1) to the set of permutations
$\sigma\in S_m\times S_{n-1}$. Further, it can be
easily checked that setting $y_n=-a_n$
in $f_{\lambda}(x/y\ts|a)$, we get
$$
f_{\lambda}(x/y\ts|a)|_{y_n=-a_n}=f_{\lambda}(x/y'\ts|a)\ts
(y_1+a_n)\cdots (y_{n-1}+a_n).
$$
On the other hand,
$$
\Delta(y)|_{y_n=-a_n}=\Delta(y')\ts
(y_1+a_n)\cdots (y_{n-1}+a_n).
$$
The factor $(y_1+a_n)\cdots (y_{n-1}+a_n)$ is symmetric in $y'$,
so removing it in the numerator and denominator we see that the
result is $\varphi_{\lambda}(x/y'\ts|a)$ which proves (5.3).

Thus, the properties (4.5), (4.6) and (4.11) are
satisfied by the polynomials $\varphi_{\lambda}(x/y\ts|a)$. So,
repeating the argument that was used for the proof of statement (i)
of Theorem 4.4, we obtain that these polynomials also satisfy (4.13).
Hence, for any partition $\zeta\subset(m,n)$-hook such that
$|\zeta|\leq |\lambda|$ and $\zeta\ne\lambda$ we have
$$
s_{\lambda}(a_{\xi}/a^*_{\eta}\ts|a)=
\varphi_{\lambda}(a_{\xi}/a^*_{\eta}\ts|a)=0.
$$
Therefore, to apply Theorem 4.5 to the polynomials
$s_{\lambda}(x/y\ts|a)$ and $\varphi_{\lambda}(x/y\ts|a)$ it remains to
check
that
$$
s_{\lambda}(a_{\xi}/a^*_{\eta}\ts|a)=
\varphi_{\lambda}(a_{\xi}/a^*_{\eta}\ts|a)
$$
for $\zeta=\lambda$. In this case $\eta=\nu$ and $\xi=\rho+\mu$.
Due to the specialization properties
(4.11) and (5.3) we may assume that $l(\nu)=n$, which implies
that the partition $\rho$ coincides with $(n^m)$. In this case we
clearly have (see also Corollary 5.2 below)
$$
\varphi_{\lambda}(x/y\ts|a)=s_{\nu}(y|a^*)s_{\mu}(x|\tau^n a)
\prod_{i=1}^m\prod_{j=1}^n(x_i+y_j).
$$
Setting $x=a_{\xi}=(\tau^na)_{\mu}$ and $y=a^*_{\nu}$ we see
that the result coincides with (4.17) which completes the proof
of the theorem.
\medskip

\noindent
{\bf Remark.} As it was noticed in [17], formula (5.1) implies that
the polynomials $s_{\lambda}(x/y\ts |a)$ coincide
with the {\it multi-Schur functions} (see, e.g., [22]) in
appropriate variables.
\medskip

As a corollary of Theorem 5.1 we obtain
the following factorization theorem
for the polynomials $s_{\lambda}(x/y\ts|a)$
which turns into
the Berele--Regev formula
for $a=(0)$ [2] (cf. [9, 35]).
\bigskip

\proclaim
{\bf Corollary 5.2} If a partition $\lambda$ is contained in the
$(m,n)$-hook
and contains the partition $(n^m)$ then
$$
s_{\lambda}(x/y\ts |a)=s_{\mu}(x|\tau^n a)s_{\nu}(y|a^*)
\prod_{i=1}^m\prod_{j=1}^n(x_i+y_j).\tag 5.4
$$
\endproclaim

\Proof Note that $\rho=(n^m)$ and that the product
$$
\prod_{(i,j)\in\rho}(x_i+y_j)=\prod_{i=1}^m\prod_{j=1}^n(x_i+y_j)
$$
is symmetric in $x$ and $y$. So, the assertion follows from (0.1).
\medskip

In the special case $\lambda=(n^m)$ formula (5.4) turns into
$$
s_{(n^m)}(x/y\ts |a)=\prod_{i=1}^m\prod_{j=1}^n(x_i+y_j). \tag 5.5
$$
Using definition (1.4) we derive from (5.5) the following analogue of the
dual Cauchy formula which is proved in [21, formula (6.17)] and 
can be also deduced from the Cauchy formula for the double
Schubert polynomials [17, 18].
\bigskip

\proclaim
{\bf Corollary 5.3}
$$
\prod_{i=1}^m\prod_{j=1}^n(x_i+y_j)=
\sum_{\lambda}s_{\tilde{\lambda}}(x|a)s_{\lambda'}(y|{-a}),\tag 5.6
$$
summed over all partitions $\lambda\subset(n^m)$, where
$\tilde{\lambda}=(n-\lambda_m,\dots,n-\lambda_1)$ and $a$ is an arbitrary
sequence.
\endproclaim

\Proof By (1.4),
$$
s_{(n^m)}(x/y\ts |{-a^*})=\sum_{\lambda\subset(n^m)}
s_{(n^m)/\lambda}(x|{-a^*})s_{\lambda'}(y|{-a}).\tag 5.7
$$
Since $s_{(n^m)/\lambda}(x|{-a^*})$ is symmetric in $x$, replacing
$x$ with $\wt x=(x_m,\dots,x_1)$ and
using (1.1) we may write
$$
s_{(n^m)/\lambda}(x|{-a^*})=s_{(n^m)/\lambda}(\wt x|{-a^*})
$$
$$
=\sum_T\prod_{\alpha\in(n^m)/\lambda}
(x_{m-T(\alpha)+1}-a_{n-T(\alpha)-c(\alpha)+1}),\tag 5.8
$$
summed over semistandard $(n^m)/\lambda$-tableaux $T$.
Consider the bijection between the cells of the diagram $(n^m)/\lambda$ and
the cells of the diagram $\tilde{\lambda}$ such that
$\alpha=(i,j)\in(n^m)/\lambda$ corresponds to
$\beta=(m-i+1,n-j+1)\in\tilde{\lambda}$. Obviously,
the content $c(\beta)$ of the cell
$\beta\in\tilde{\lambda}$ is related with $c(\alpha)$ by
$c(\beta)=n-m-c(\alpha)$.
Moreover, the map
$$
T(\alpha)\to \wt T(\beta)=m-T(\alpha)+1
$$
is a bijection between the semistandard $(n^m)/\lambda$-tableaux and
the semistandard $\tilde{\lambda}$-tableaux. So, (5.8) gives
$$
s_{(n^m)/\lambda}(x|{-a^*})=
\sum_{\wt T}\prod_{\beta\in\tilde{\lambda}}
(x_{\wt T(\beta)}-a_{\wt T(\beta)+
c(\beta)})=s_{\tilde{\lambda}}(x|a).
$$
Using (5.7) and (5.5) we complete the proof.

\bigskip
\bigskip
\noindent
{\bf 6. Macdonald--Goulden--Greene formula}
\bigskip

Formula (5.5) means that the polynomial $s_{(n^m)}(x/y\ts |a)$
does not depend on the sequence $a$. It turns out, that an analogous
phenomenon takes place for infinite families
of variables. We shall regard the elements of the sequence $a$
as independent variables to avoid convergence problems.
Let us consider three families of variables $x=(x_i)$,
$y=(y_i)$, $a=(a_i)$, $i\in\ZZ$. We define
the functions $s_{\lambda/\mu}(x/y\ts |a)$
in $x$, $y$ and $a$
by formula (1.5) where we allow the primed and unprimed
entries of the
tableaux to run through the set of all integers. 
This definition
can be shown to be
equivalent to the following formula, where the factorial Schur
functions $s_{\lambda/\mu}(x|a)$ are defined by (1.1) with
$T$ running over semistandard $\lambda/\mu$-tableaux with
entries from $\ZZ$ (see [9] and [21]).
\bigskip

\proclaim
{\bf Proposition 6.1}
$$
s_{\lambda/\mu}(x/y\ts|a)=\sum_{\mu\subset\nu\subset\lambda}
s_{\lambda/\nu}(x|a)s_{\nu'/\mu'}(y|{-a}).\tag 6.1
$$
\endproclaim

\Proof Repeating the arguments of the proof of
Proposition 1.2 we obtain that the assertion follows from the formula
$$
\sum_{T}\prod_{\alpha\in\nu/\mu}
(y_{T(\alpha)}+a_{T(\alpha)+c(\alpha)})=s_{\nu'/\mu'}(y|{-a}),\tag 6.2
$$
summed over $\nu/\mu$-tableaux $T$ with entries from $\ZZ$
whose rows strictly decrease and
columns weakly decrease.

It was shown in [9] and [21] that in the case of infinite number
of variables (parametrized by $\ZZ$)
the factorial and supersymmetric Schur functions coincide with each other:
$$
s_{\lambda/\mu}(x|a)=s_{\lambda/\mu}(x/{-a}),\tag 6.3
$$
where the supersymmetric Schur functions $s_{\lambda/\mu}(x/y)$ are
still defined by (0.2). 
In particular, $s_{\lambda/\mu}(x|a)$ is symmetric in $a$ (which is not
true
in the finite case).
Hence, we have 
$$
s_{\nu'/\mu'}(y|{-a})=s_{\nu'/\mu'}(\wt y\ts|{-\wt a}),
$$
where $\wt y=(y_{-i})$ and $\wt a=(a_{-i})$.
So, by (1.1),
$$
s_{\nu'/\mu'}(y|{-a})
=\sum_{T'}\prod_{\alpha'\in\nu'/\mu'}
(y_{-T'(\alpha')}+a_{-T'(\alpha')-c(\alpha')}),\tag 6.4
$$
summed over semistandard $\nu'/\mu'$-tableaux $T'$ with entries from
$\ZZ$. Note that the map
$$
T'(\alpha')\to T(\alpha)=-T'(\alpha'),
$$
where $\alpha=(i,j)\in\nu/\mu$ and $\alpha'=(j,i)\in\nu'/\mu'$, is
a bijection between the set of semistandard $\nu'/\mu'$-tableaux
and the set of $\nu/\mu$-tableaux whose rows strictly decrease and
columns weakly decrease. Obviously, $c(\alpha)=-c(\alpha')$, hence,
(6.4) coincides with the left hand side of (6.2) which completes the
proof.
\medskip

For $a=(0)$ formula (6.1) turns into the definition of
the supersymmetric Schur functions $s_{\lambda/\mu}(x/y)$. 
It turns out that the right hand side of (6.1) does not depend on
the variables $a_i$ and thus the Macdonald--Goulden--Greene formula
(see [9] and [21])
holds for the functions  $s_{\lambda/\mu}(x/y|a)$ as well.

\proclaim
{\bf Theorem 6.2} One has the formula
$$
s_{\lambda/\mu}(x/y|a)=\sum_T\prod_{\alpha\in\lambda/\mu}
(x_{T(\alpha)}+y_{T(\alpha)+c(\alpha)}),\tag 6.5
$$
where $T$ runs over all semistandard $\lambda/\mu$-tableaux with entries
from $\ZZ$.
\endproclaim

\Proof For $a=(0)$ formula (6.5) was proved in [9] and [21].
So, it suffices to
verify that $s_{\lambda/\mu}(x/y|a)=s_{\lambda/\mu}(x/y)$. Using (6.1) and
(6.3) we obtain
$$
s_{\lambda/\mu}(x/y|a)=\sum_{\nu}s_{\lambda/\nu}(x/{-a})
s_{\nu'/\mu'}(y/a).\tag 6.6
$$
By the symmetry property (0.3) we have $s_{\nu'/\mu'}(y/a)=
s_{\nu/\mu}(a/y)$. Hence, using
the definition of the supersymmetric Schur functions we can rewrite (6.6)
as follows
$$
\align
s_{\lambda/\mu}(x/y|a)&=\sum_{\nu,\rho,\sigma}s_{\lambda/\rho}(x)
s_{\rho'/\nu'}({-a})s_{\nu/\sigma}(a)s_{\sigma'/\mu'}(y)\\
&=\sum_{\rho,\sigma}s_{\lambda/\rho}(x)s_{\sigma'/\mu'}(y)
s_{\rho'/\sigma'}({-a}/a).\tag 6.7
\endalign
$$
Note now that $s_{\rho'/\sigma'}({-a}/a)=0$ unless $\rho=\sigma$. Indeed,
suppose that there exists a cell $\alpha_0\in\rho'/\sigma'$. Then,
applying (6.5) with $a=(0)$ we get
$$
s_{\rho'/\sigma'}({-a}/a)=s_{\rho'/\sigma'}({-a}/\tau^{-c(\alpha_0)}a)
$$
$$
=\sum_T\prod_{\alpha\in\rho'/\sigma'}
({-a}_{T(\alpha)}+a_{T(\alpha)+c(\alpha)-c(\alpha_0)})=0.
$$
Thus, (6.7) coincides with $s_{\lambda/\mu}(x/y)$ which completes the
proof.

\bigskip
\bigskip
\noindent
{\bf 7. Shifted supersymmetric Schur
polynomials and a basis in the center of $\U(\gl(m|n))$}
\bigskip

We need to introduce super-analogues of the shifted symmetric polynomials
(cf. [29, 31, 32]). A polynomial in two families of variables
$u=(u_1,\dots,u_m)$ and $v=(v_1,\dots,v_n)$ will be called
{\it shifted supersymmetric} if it is supersymmetric in
the variables
$$
(u_1+m-1,\ts u_2+m-2,\ts\dots,\ts u_m)
\qquad
\text{and}
\qquad
(v_1,\ts v_2-1,\ts\dots,\ts v_n-n+1).
$$
We shall
denote the algebra of shifted supersymmetric polynomials by
$\Lambda^*(m|n)$.
It follows from its definition that $\Lambda^*(m|n)$
is isomorphic to the algebra of supersymmetric
polynomials in $x$ and $y$.

Let us consider the polynomials $s_{\lambda/\mu}(x/y\ts |a)$
with the
sequence $a$ defined by $a_i=-m+i$ and
put in $s_{\lambda/\mu}(x/y\ts |a)$ 
$$
\alignat2
x_i&=u_{m-i+1}-m+i\qquad&&\text{for}\quad i=1,\dots,m,\\
y_j&=v_j+m-j\qquad&&\text{for}\quad j=1,\dots,n.
\endalignat
$$
Then we obtain a shifted supersymmetric polynomial in $u$ and $v$
which will be denoted by $s^*_{\lambda/\mu}(u/v)$ and will be called
{\it shifted supersymmetric Schur polynomial}. 
In the case $n=0$
it coincides with the shifted Schur polynomial $s^*_{\lambda/\mu}(u)$
(see [29, 31, 32]) which can be defined by the formula
$$
s^*_{\lambda/\mu}(u)=\sum_{T}
\prod_{\alpha\in\lambda/\mu}
(u_{T(\alpha)}-c(\alpha)),\tag 7.1
$$
summed over $\lambda/\mu$-tableaux $T$ with entries in
$\{1,\dots,m\}$ whose rows
weakly decrease and columns strictly decrease.
The highest component of $s^*_{\lambda/\mu}(u)$ is the usual skew Schur
polynomial $s_{\lambda/\mu}(u)$.

We shall state now some properties of the polynomials
$s^*_{\lambda/\mu}(u/v)$,
which can be easily derived from the corresponding properties of the
polynomials $s_{\lambda/\mu}(x/y\ts |a)$.
\medskip

First we give a combinatorial interpretation of $s^*_{\lambda/\mu}(u/v)$.

To distinguish the indices of $u$ and $v$ let us identify
the indices of $v$ with the symbols $1',\dots,n'$.
Consider the diagram of shape $\lambda/\mu$ and fill it with the
indices $1',\dots,n',1,\dots,m$ such that:
\medskip
(a) In each row (resp. column) each primed index is to the left (resp.
above)
from each unprimed index.

(b) Primed indices strictly decrease along rows and weakly decrease
down columns.

(c) Unprimed indices weakly decrease along rows and strictly decrease
down columns.
\medskip

Denote the resulting tableau by $T$.

\proclaim
{\bf Proposition 7.1} One has the formula
$$
s^*_{\lambda/\mu}(u/v)=\sum_{T}
\prod_{\underset{\ssize T(\alpha)\ts \text{unprimed}}
\to{\alpha\in\lambda/\mu}}
(u_{T(\alpha)}-c(\alpha))
\prod_{\underset{\ssize T(\alpha)\ts \text{primed}}
\to{\alpha\in\lambda/\mu}}
(v_{T(\alpha)}+c(\alpha)).\tag 7.2
$$
\endproclaim

\Proof This follows immediately from Proposition 1.2. It suffices
to use (1.5) with $x$ replaced by $\wt x=(x_m,\dots,x_1)$.
\medskip

In particular, for the {\it elementary\/} and {\it complete
shifted supersymmetric polynomials} \newline
$e^*_k(u/v):=s^*_{(1^k)}(u/v)$ and
$h^*_k(u/v):=s^*_{(k)}(u/v)$ we have
$$
e^*_k(u/v)=\sum_{p+q=k}
\sum_{\underset{\ssize j_1\geq\cdots\geq j_q}\to
{i_1>\cdots >i_p}}v_{j_1}(v_{j_2}-1)\cdots (v_{j_q}-q+1)
(u_{i_1}+q)\cdots (u_{i_p}+k-1),
$$
$$
h^*_k(u/v)=\sum_{p+q=k}
\sum_{\underset{\ssize j_1>\cdots > j_q}\to
{i_1\geq\cdots \geq i_p}}v_{j_1}(v_{j_2}+1)\cdots (v_{j_q}+q-1)
(u_{i_1}-q)\cdots (u_{i_p}-k+1).
$$
Using (7.1) we can rewrite (7.2) in the following equivalent form.

\bigskip
\proclaim
{\bf Corollary 7.2} One has
$$
s^*_{\lambda/\mu}(u/v)=\sum_{\mu\subset\nu\subset\lambda}
s^*_{\lambda/\nu}(u)\ts s^*_{\nu'/\mu'}(v).\tag 7.3
$$
\endproclaim

This implies that the highest component of $s^*_{\lambda/\mu}(u/v)$
is the supersymmetric Schur polynomial $s_{\lambda/\mu}(u/v)$.
So, the polynomials $s^*_{\lambda}(u/v)$ with
$\lambda\subset(m,n)$-hook form a basis in $\Lambda^*(m|n)$.

The following are reformulations of Theorems 4.4 and 4.5$'$
for the shifted supersymmetric polynomials; cf. [29, 31]
(we use the notation from Section 4).
\bigskip

\proclaim
{\bf Theorem 7.3} Let $\lambda,\zeta$ be partitions which are contained in
the $(m,n)$-hook. 

{\rm (i)} If $\lambda\not\subset\zeta$ then
$$
s^*_{\lambda}(\xi/\eta)=0. \tag 7.4
$$

{\rm (ii)} If $\lambda=\zeta$ then
$$
s^*_{\lambda}(\xi/\eta)=H(\lambda), \tag 7.5
$$
where $H(\lambda)$ is the product of the hook lengths of all cells of
$\lambda$.
\endproclaim

\proclaim
{\bf Theorem 7.4} Let $f(u/v)$ be a shifted
supersymmetric polynomial such that
$$
f(u/v)=s_{\lambda}(u/v)+\text{lower terms}
$$
for some partition $\lambda\subset(m,n)$-hook, and
$
f(\xi/\eta)=0
$
for any partition $\zeta\subset(m,n)$-hook with $|\zeta|<|\lambda|$. Then
$f(u/v)=s^*_{\lambda}(u/v)$.
\endproclaim
\bigskip

A distinguished linear basis in the center of the universal enveloping
algebra
$\U(\gl(m))$ was constructed in [29]. The eigenvalue of a basis element
in a highest weight representation is a shifted Schur polynomial
$s^*_{\lambda}(u)$. It turns out, that this construction can be easily
carried to the case of the Lie superalgebra $\gl(m|n)$.
Below we formulate the corresponding theorem and briefly outline
its proof. Another approach to this construction, based on the
super-analogues of the higher
Capelli identities is contained in Section 8.
\medskip

We denote by $E_{ij}$, $i,j=1,\dots,m+n$ the standard basis of
the Lie superalgebra $\gl(m|n)$. The $\ZZ_2$-grading on $\gl(m|n)$
is defined by $E_{ij}\mapsto p(i)+p(j)$, where $p(i)=0$ or $1$
depending on whether $i\leq m$ or $i>m$. 
The commutation relations in this basis
are given by
$$
[E_{ij},E_{kl}]=\delta_{kj}E_{il}-\delta_{il}E_{kj}
(-1)^{(p(i)+p(j))(p(k)+p(l))}.\tag 7.6
$$

Given $w=(u_1,\dots,u_m,v_1,\dots,v_n)\in\C^{m+n}$ we consider an
arbitrary highest weight $\gl(m|n)$-module $L(w)$
with the highest weight $w$.
That is, $L(w)$ is generated by a nonzero vector $\psi$ such that
$$
\alignat2
E_{ii}\ts\psi&=u_i\ts\psi\qquad&&\text{for}\quad i=1,\dots,m,\\
E_{m+j,m+j}\ts\psi&=v_j\ts\psi\qquad&&\text{for}\quad j=1,\dots,n,\\
E_{ij}\ts\psi&=0\qquad&&\text{for}\quad 1\leq i<j\leq m+n.
\endalignat
$$
Every element $z$ of the center
$\Z(\gl(m|n))$ of the universal enveloping algebra
$\U(\gl(m|n))$ acts in $L(w)$ as a scalar $\chi(z)$.
For a fixed $z$ the scalar $\chi(z)$ is a shifted supersymmetric
polynomial in $u$ and $v$ and
the map
$z\mapsto \chi(z)$ defines an algebra isomorphism
$$
\chi:\Z(\gl(m|n))\to \Lambda^*(m|n),\tag 7.7
$$
which is called the
Harish-Chandra isomorphism; see [16, 38, 40].

Our aim now is to give an explicit description of the basis of
the algebra $\Z(\gl(m|n))$ formed by the preimages
$\chi^{-1}(s^*_{\lambda}(u/v))$ of the basis elements of $\Lambda^*(m|n)$.
Let us introduce some more notations.

We shall need to consider matrices with entries from superalgebras.
All our matrices will be even. That is, if
$B=(B_{ia})$ is a $(m,n)\times(m',n')$-matrix whose entries are homogeneous
entries of a superalgebra $\Cal B$, 
we shall always have $p(B_{ia})=p(i)+p(a)$, where
$p(a)=0$ or $1$ depending on whether $a\leq m'$ or $a>m'$. A matrix $B$
will
be identified with an element of the tensor product
$$
B=\sum_{i,a}e_{ia}\ot B_{ia}(-1)^{p(a)(p(i)+1)}
\in\Mat_{(m,n)\times (m',n')}\ot\Cal B,
$$
where the $e_{ia}$ are the standard matrix units.

More generally,
given $k$ matrices $B^{(1)},\dots, B^{(k)}$ of the size
$(m,n)\times(m',n')$
we define their tensor product $B^{(1)}\ot\cdots\ot B^{(k)}$
as an element
$$
\sum e_{i_1a_1}\ot \cdots\ot e_{i_ka_k}\ot B^{(1)}_{i_1a_1}\cdots 
B^{(k)}_{i_ka_k}\ts (-1)^{\gamma(I,A)}
\in(\Mat_{(m,n)\times (m',n')})^{\ot k}\ot\Cal B,
$$
where
$$
\gamma(I,A)=
\sum_{r=1}^k p(a_r)(p(i_r)+1)+\sum_{1\leq r<s\leq k}
(p(i_r)+p(a_r))(p(i_s)+p(a_s)).
$$

The supertrace of an element
$$
B=\sum e_{i_1j_1}\ot \cdots\ot e_{i_kj_k}\ot 
B_{i_1,\dots,i_k;j_1,\dots,j_k}\in 
(\Mat_{(m,n)\times (m,n)})^{\ot k}\ot\Cal B
$$
is defined by
$$
\str B=\sum_{i_1,\dots,i_k}B_{i_1,\dots,i_k;i_1,\dots,i_k}
(-1)^{p(i_1)+\cdots+p(i_k)}.
$$

Using the natural action of the symmetric group $S_k$ in the space
$(\C^{m|n})^{\ot k}$ we represent each element of $S_k$ by
a linear combination of tensor products of matrices.
In particular, the transposition $(i,j)\in S_k$, $i<j$,
corresponds to the element
$$
P_{ij}=\sum_{a,b} 1\ot\cdots\ot 1\ot e_{ab}\ot 1\ot\cdots\ot 1\ot
e_{ba}\ot 1\ot\cdots\ot 1\ts (-1)^{p(b)},
$$
where the tensor factors $e_{ab}$ and $e_{ba}$ take the $i$th and $j$th
places, respectively.

We can now describe the construction of a basis in $\Z(\gl(m|n))$.

Set $\wh E_{ij}=E_{ij}(-1)^{p(j)}$ and denote by $\wh E$ the
$(m,n)\times (m,n)$-matrix
whose $ij$th entry is $\wh E_{ij}$.

Following [29], for a partition
$\lambda$ and a standard $\lambda$-tableau $T$ we denote by $v_T$
the corresponding vector of the Young orthonormal basis
with respect to an invariant inner product $(\ ,\ )$ in the
irreducible representation $V^{\lambda}$ of the symmetric group $S_k$,
$k=|\lambda|$. We let $c_T(r)=j-i$ if the cell $(i,j)\in\lambda$
is occupied by the entry $r$ of the tableau $T$. Given two
standard $\lambda$-tableaux $T$ and $T'$,
introduce the matrix element
$$
\Psi_{TT'}=\sum_{s\in S_k}(s\cdot v_T,v_{T'})\cdot s^{-1}\in\C[S_k].
\tag 7.8
$$

\bigskip
\proclaim
{\bf Theorem 7.5} The element
$$
\SS_{\lambda}=\frac{1}{H(\lambda)}\ts\str 
(\wh E-c_T(1))\ot\cdots\ot(\wh E-c_T(k))\cdot\Psi_{TT},\tag 7.9
$$
is independent of a $\lambda$-tableau $T$. The set of elements
$\SS_{\lambda}$
with $\lambda\subset(m,n)$-hook forms a basis in $\Z(\gl(m|n))$.
Moreover, the image of $\SS_{\lambda}$ under the Harish-Chandra
isomorphism is $s^*_{\lambda}(u/v)$.
\endproclaim

\noindent
{\bf Outline of the proof.}
In the case $n=0$ this theorem was proved in [29]. It
constitutes the `difficult part' of the proof of the higher
Capelli identities.
A straightforward generalization of those
arguments proves Theorem 7.5.
For this one
uses the following
$R$-matrix form of the defining relations in $\U(\gl(m|n))$ (cf. [27]):
$$
R(u-v)\cdot \wh E(u)\ot\wh E(v)=\wh E(v)\ot\wh E(u)\cdot R(u-v),
$$
where
$
R(u)=1+P_{12}\ts u.
$

\bigskip
\noindent
{\bf Examples.} For the partitions of weight $\leq 2$ we have
$$
\SS_{(1)}=\str \wh E=\sum_i E_{ii},
$$
$$
\align
\SS_{(2)}&=\frac12\ts \str\bigl(\wh E\ot(\wh E-1)\cdot (1+P_1)\bigr)\\
&=\frac12\ts \sum_{i,j}\bigl(E_{ii}(E_{jj}-(-1)^{p(j)})+
E_{ij}(E_{ji}-\delta_{ji}(-1)^{p(i)})(-1)^{p(j)}\bigr),
\endalign
$$
$$
\align
\SS_{(1^2)}&=\frac12\ts \str\bigl(\wh E\ot(\wh E+1)\cdot (1-P_1)\bigr)\\
&=\frac12\ts \sum_{i,j}\bigl(E_{ii}(E_{jj}+(-1)^{p(j)})-
E_{ij}(E_{ji}+\delta_{ji}(-1)^{p(i)})(-1)^{p(j)}\bigr).
\endalign
$$
\bigskip
\bigskip
\noindent
{\bf 8. Super Capelli identities}
\bigskip

Here we formulate a super-analogue of the higher Capelli identities
obtained in [28--30].

Consider the supercommutative algebra $\Cal Z$ with generators $z_{ia}$
where
$i=1,\dots,m+n$ and $a=1,\dots,m'+n'$ and the $\ZZ_2$-grading
given by $z_{ia}\mapsto p(i)+p(a)$. Define the representation $\pi$
of the Lie superalgebra $\gl(m|n)$ in $\Cal Z$ by
$$
\pi(E_{ij})=\sum_{a=1}^{m'+n'}z_{ia}\di_{ja},
$$
where $\di_{ja}=\di/\di z_{ja}$ is the left derivation.
This definition can be rewritten in the matrix form
as follows:
$$
\pi(\wh E)=Z\ts D',\tag 8.1
$$
where $Z$ is the $(m,n)\times(m',n')$-matrix $(z_{ia})$ and
and $D'$ is the $(m',n')\times(m,n)$-matrix $(\di'_{ai})$
with $\di'_{ai}=\di_{ia}(-1)^{p(i)}$.

For a partition $\lambda$ with $|\lambda|=k$ denote by $\chi^{\lambda}$
the irreducible character of $S_k$. We identify $\chi^{\lambda}$
with an element of the group algebra $\C[S_k]$:
$$
\chi^{\lambda}=\sum_{s\in S_k}\chi^{\lambda}(s)\cdot s.
$$
Define the differential operator $\Delta_{\lambda}$ by
$$
\Delta_{\lambda}=\frac 1{k!}\ts
\str\bigl(Z^{\ot k}\cdot {D'}^{\ot k}\cdot \chi^{\lambda}\bigr).
$$

The following is a super-analogue of the higher Capelli identities
(cf. [28--30]).
\bigskip
\proclaim
{\bf Theorem 8.1}
$$
\pi(\SS_{\lambda})=\Delta_{\lambda}.\tag 8.2
$$
\endproclaim

We outline two proofs of this identity. The first proof is based on
the properties of the shifted supersymmetric polynomials.
The second one uses a super-analogue of a more general identity
obtained in [28] and [30] (see Theorem 8.2 below).
\medskip

\noindent
{\bf The first proof.} We follow again the corresponding arguments
from [29].
First, one
verifies that the operator $\Delta_{\lambda}$ commutes with
both the actions of the Lie superalgebras $\gl(m|n)$ and $\gl(m'|n')$
in $\Cal Z$; the latter is given by
$$
\pi'(E'_{ab})=\sum_{i=1}^{m+n}z_{ia}\di_{ib}(-1)^{(p(a)+p(b))p(i)},
$$
where the $E'_{ab}$ are the standard generators of $\gl(m'|n')$.
This implies that $\Delta_{\lambda}$ is the image of
a certain element $\SS'_{\lambda}\in\Z(\gl(m|n))$ under $\pi$
(cf. [12, 27]). To prove that
$\SS_{\lambda}=\SS'_{\lambda}$ we compare their images under
the Harish-Chandra isomorphism. Set
$s'_{\lambda}(u/v)=\chi(\SS'_{\lambda})$.
By Theorem 7.5, $\chi(\SS_{\lambda})=s^*_{\lambda}(u/v)$.
We use Theorem 7.4 to prove that $s'_{\lambda}(u/v)=s^*_{\lambda}(u/v)$.
Using (8.1)
we check that both sides of (8.2) agree modulo lower terms (cf. [29, 31]).
This proves that the polynomials
$s'_{\lambda}(u/v)$ and $s^*_{\lambda}(u/v)$ have
the same highest component which coincides with
the supersymmetric Schur polynomial $s_{\lambda}(u/v)$.

By Theorem 7.3, to complete the proof we have to verify that
$s'_{\lambda}(\xi/\eta)$ is zero for any partition $\zeta\subset(m,n)$-hook
such that $|\zeta|<|\lambda|$. 
Let us consider the superalgebra $\Cal Z$ with the parameters $m'$ and $n'$
being sufficiently large, so that $m'\geq\text{max}\{m,|\zeta|\}$ and
$n'\geq\text{max}\{n,|\zeta|\}$.

Introduce the
following element of $\Cal Z$:
$$
\psi_{\zeta}=\Delta_1^{\eta_1-\eta_2}
\Delta_2^{\eta_2-\eta_3}\cdots
\Delta_n^{\eta_n}
\ts\prod_{(i,j)\in\xi}z_{i,m'+j},
$$
where ${\Delta}_r=\det[z_{m+i,m'+j}]_{1\leq i,j\leq r}$ and
the product is taken in any fixed order.
It can be easily checked that $\psi_{\zeta}$ satisfies the
relations
$$
\alignat2
\pi(E_{ii})\ts \psi_{\zeta}&=\xi_i\ts \psi_{\zeta}\qquad&&\text{for}\quad
i=1,\dots,m,\\
\pi(E_{m+j,m+j})\ts \psi_{\zeta}&=\eta_j\ts
\psi_{\zeta}\qquad&&\text{for}\quad
j=1,\dots,n,\\
\pi(E_{ij})\ts \psi_{\zeta}&=0\qquad&&\text{for}\quad
1\leq i<j\leq m+n.
\endalignat
$$
This means that $\psi_{\zeta}$ generates a 
$\gl(m|n)$-module with
the highest weight $(\xi,\eta)$. Hence, $\psi_{\zeta}$
is an eigenvector for the operator $\Delta_{\lambda}$ with
the eigenvalue $s'_{\lambda}(\xi/\eta)$. However, the degree of
$\psi_{\zeta}$ equals $|\zeta|$ and so, if $|\zeta|<|\lambda|=k$ then
$\psi_{\zeta}$ is annihilated by $\Delta_{\lambda}$, that is,
$s'_{\lambda}(\xi/\eta)=0$ which completes the proof.

\bigskip
\proclaim
{\bf Theorem 8.2} Let $T$ and $T'$ be
two standard tableaux of the shape $\lambda$. Then
$$
\pi\left((\wh E-c_T(1))\ot\cdots
\ot(\wh E-c_T(k))\cdot\Psi_{TT'}\right)=
Z^{\ot k}\cdot(D')^{\ot k} \cdot\Psi_{TT'}.\tag 8.3
$$
\endproclaim

\Proof Following [30], we use some properties of the Jucys--Murphy
elements in the group algebra for the symmetric group.
However, the
arguments from [30] can be modified to avoid using
the Wick formula and the Olshanski\u\i\
special symmetrization map.
\medskip

We use induction on $k$. Denote by $U$ the tableau
obtained from $T$ by removing the cell with the entry $k$.
From the branching property of the Young basis $\{v_T\}$ one
can easily derive that
$$
\Psi_{TT'}=\text{\rm const}\cdot\Psi_{UU}\Psi_{TT'},
$$
where `const' is a nonzero constant (more precisely, $\text{const}=
\text{dim}\ts \mu/(k-1)!$ where $\mu$ is the shape of $U$
and $\text{dim}\ts \mu=\text{dim}\ts V^{\mu}$).

So, we can rewrite the
left hand side of (8.3) as follows:
$$
\text{const}\cdot(ZD'-c_T(1))\ot\cdots\ot(ZD'-c_T(k-1))\cdot\Psi_{UU}
\ot(ZD'-c_T(k))\cdot\Psi_{TT'}.
$$
By the induction hypothesis, this equals
$$
\text{\rm const}\cdot Z^{\ot k-1}\cdot(D')^{\ot k-1} \cdot\Psi_{UU}
\ot(ZD'-c_T(k))\cdot\Psi_{TT'}
$$
$$
=Z^{\ot k-1}\cdot(D')^{\ot k-1}
\ot(ZD'-c_T(k))\cdot\Psi_{TT'}
$$
$$
=\bigl(\sum e_{i_1j_1}\ot\cdots\ot e_{i_kj_k}\ot
z_{i_1a_1}\cdots z_{i_{k-1}a_{k-1}}
\di'_{a_1j_1}\cdots \di'_{a_{k-1}j_{k-1}}
$$
$$
(\sum z_{i_ka_k}\di'_{a_kj_k}-\delta_{i_kj_k}c_T(k))
\bigr)(-1)^{\alpha(I,J,A)}\cdot\Psi_{TT'},
$$
where
$$
\align
\alpha(I,J,A)=\sum_{r=1}^k p(j_r)(p(i_r)+1)+&\sum_{1\leq r<s\leq k-1}
(p(a_r)+p(j_r))(p(i_s)+p(a_s))\\
+&\sum_{1\leq r<s\leq k}
(p(i_r)+p(j_r))(p(i_s)+p(j_s)).
\endalign
$$
Now we transform this expression using the relations 
$$
\di'_{bj}z_{ia}=z_{ia}\di'_{bj}(-1)^{(p(i)+p(a))(p(j)+p(b))}+
\delta_{ab}\delta_{ij}(-1)^{p(j)}
$$
to obtain
$$
\left(\sum e_{i_1j_1}\ot\cdots\ot e_{i_kj_k}\ot
z_{i_1a_1}\cdots z_{i_{k}a_{k}}
\di'_{a_1j_1}\cdots \di'_{a_{k}j_{k}}(-1)^{\beta_k(I,J,A)}
\right)\cdot\Psi_{TT'}
$$
$$
+\left(\sum e_{i_1j_1}\ot\cdots\ot e_{i_{k-1}j_{k-1}}\ot 1\ot
z_{i_1a_1}\cdots z_{i_{k-1}a_{k-1}}
\di'_{a_1j_1}\cdots \di'_{a_{k-1}j_{k-1}}(-1)^{\beta_{k-1}(I,J,A)}
\right)
$$
$$
\times(P_{1k}+\cdots+P_{k-1,k}-c_T(k))\cdot\Psi_{TT'},\tag 8.4
$$
where
$$
\align
\beta_k(I,J,A)=\sum_{r=1}^k p(j_r)(p(i_r)+1)+&\sum_{1\leq r<s\leq k}
(p(a_r)+p(j_r))(p(i_s)+p(a_s))\\
+&\sum_{1\leq r<s\leq k}
(p(i_r)+p(j_r))(p(i_s)+p(j_s)).
\endalign
$$
Note that $P_{1k}+\cdots+P_{k-1,k}$ is the image of the Jucys--Murphy
element $(1k)+\cdots+(k-1,k)\in\C[S_k]$ (see [14] and
[26]). It has the property
$$
((1k)+\cdots+(k-1,k))\cdot\Psi_{TT'}=c_T(k)\cdot\Psi_{TT'},
$$
which was also used in [30] and can be easily derived from the following
formula due to Jucys and Murphy:
$$
((1k)+\cdots+(k-1,k))\cdot v_T=c_T(k)\cdot v_T.
$$
This proves that the second summand in (8.4)
is zero, while the first one coincides with the right hand side
of (8.3).
Theorem 8.2 is proved.
\medskip

\noindent
{\bf The second proof of Theorem 8.1.} Put $T=T'$ in Theorem 8.2 and take
the supertrace of both sides of (8.3). On the left hand side we
get $H(\lambda)\ts \pi(\SS_{\lambda})$ while for the right hand side
we have
$$
\align
\str Z^{\ot k}\cdot(D')^{\ot k} \cdot\Psi_{TT}&=
\frac{1}{k!}\sum_{s\in S_k}\str s\cdot Z^{\ot k}\cdot(D')^{\ot k} 
\cdot\Psi_{TT}\cdot s^{-1}\\
&=\frac{1}{\text{dim}\ts \lambda}\ts\str Z^{\ot k}\cdot(D')^{\ot k} \cdot
\chi^{\lambda}.
\endalign
$$
Here we have used the invariance of $Z^{\ot k}$ and $(D')^{\ot k}$
under the conjugations by elements $s\in S_k$ and the following
equality of elements of the group algebra of $S_k$:
$$
\frac{1}{k!}\sum_{s\in S_k} s\cdot \Psi_{TT}\cdot s^{-1}=
\frac{1}{\text{dim}\ts \lambda}\ts\chi^{\lambda}.
$$
So, on the right hand side we get $H(\lambda)\ts \Delta_{\lambda}$
which completes the proof.

\bigskip
\noindent
{\bf Example.} In the case of $\lambda=(1^k)$ the identity (8.2) was
obtained by M. Nazarov [27] in another form. Namely,
a formal series $B(t)$ whose coefficients are generators of $\Z(\gl(m|n))$
was constructed in [27] with the use of some
properties of the Yangian for the Lie superalgebra $\gl(m|n)$.
The explicit expression for $B(t)$ has the form
of a `quantum' analogue of Berezinian:
$$
\aligned
B(t)=\sum_{\sigma\in S_m}\sgn(\sigma)&
\left(1+\frac{\wh E}{t}\right)_{\sigma(1),1}\cdots 
\left(1+\frac{\wh E}{t-m+1}\right)_{\sigma(m),m}\\
\times \sum_{\tau\in S_n}\sgn(\tau)&
\left(1+\frac{\wh E}{t-m+1}\right)^{*}_{m+\tau(1),m+1}\cdots
\left(1+\frac{\wh E}{t-m+n}\right)^{*}_{m+\tau(n),m+n},
\endaligned
\tag 8.5
$$
where $A^*=(A^{-1})^{\text{\it st}}$ and {\it st} is the matrix
supertransposition:
$(B^{\text{\it st}})_{ij}=B_{ji}(-1)^{p(i)(p(j)+1)}$.
The image of $B(t)$ under the Harish-Chandra isomorphism
coincides with its eigenvalue on
the highest vector $\psi$ of the highest weight
$\gl(m|n)$-module $L(w)$, $w=(u,v)\in\C^{m|n}$. The eigenvalue of the first
determinant in (8.5) on $\psi$ is
$$
\frac{(t+u_1)\cdots (t+u_m-m+1)}
{t(t-1)\cdots (t-m+1)}.
$$
To find the eigenvalue of the second determinant we may replace the matrix
$\wh E$ with its submatrix $\wt E=(\wh E_{ij})_{m+1\leq i,j\leq m+n}$.
However, the determinant
$$
\sum_{\tau\in S_n}\sgn(\tau)
\left(1+\frac{\wt E}{t-m+1}\right)^{*}_{m+\tau(1),m+1}\cdots
\left(1+\frac{\wt E}{t-m+n}\right)^{*}_{m+\tau(n),m+n}
$$
equals
$$
\left(\sum_{\tau\in S_n}\sgn(\tau)
\left(1+\frac{\wt E}{t-m+1}\right)_{m+\tau(1),m+1}\cdots
\left(1+\frac{\wt E}{t-m+n}\right)_{m+\tau(n),m+n}\right)^{-1},
$$
which follows from [27, Proposition 3] and can be also proved directly
by using an $R$-matrix form of the defining relations in $\U(\gl(n))$;
see, e.g., [25]. So, its eigenvalue on $\psi$ is
$$
\frac{(t-m+1)\cdots (t-m+n)}
{(t-v_1-m+1)\cdots (t-v_n-m+n)}.
$$
Thus,
$$
\chi(B(t))=\frac{(t+u_1)\cdots (t+u_m-m+1)(t-m+1)\cdots (t-m+n)}
{t(t-1)\cdots (t-m+1)(t-v_1-m+1)\cdots (t-v_n-m+n)}.
$$
Relation (2.5) implies that
$$
\chi(B(t))=1+\sum_{k=1}^{\infty}\frac{e^*_k(u/v)}{t(t-1)\cdots (t-k+1)}.
$$
By Theorem 7.5, $\chi(\SS_{(1^k)})=e^*_k(u/v)$, hence,
$$
B(t)=1+\sum_{k=1}^{\infty}\frac{\SS_{(1^k)}}{t(t-1)\cdots (t-k+1)}.
$$
Using (8.2) we get the identity (see [27]):
$$
\pi(B(t))=1+\sum_{k=1}^{\infty}\frac{\Delta_{(1^k)}}{t(t-1)\cdots (t-k+1)}.
$$
For $n=n'=0$ it turns into the classical Capelli identity;
see, e.g., [11, 12].
\bigskip
\bigskip

\noindent
{\bf References}
\bigskip

\itemitem{[1]}{A. Abderrezzak}, {G\'en\'eralisation d'identit\'es de
Carlitz,
Howard et Lehmer}, {\sl Aequat. Math.} {\bf 49} (1995), 36--46.

\itemitem{[2]}{A. Berele and A. Regev}, {Hook Young diagrams with
applications to combinatorics and to representations of Lie
superalgebras}, {\sl Adv. Math.} {\bf 64} (1987), 118--175.

\itemitem{[3]}{N. Bergeron and A. M. Garsia}, {Sergeev's formula
and the Littlewood--Richardson rule}, {\sl Linear and Multilinear Algebra}
{\bf 27} (1990), 79--100.

\itemitem{[4]}
{L.~C.~Biedenharn and J.~D.~Louck},
{A new class of symmetric
polynomials defined in terms of tableaux},
{\sl Advances in Appl.\ Math.}
{\bf 10} (1989),
396--438.

\itemitem{[5]}
{L.~C.~Biedenharn and J.~D.~Louck},
{Inhomogeneous basis set of
symmetric polynomials defined by tableaux},
{\sl Proc.\ Nat.\ Acad.\ Sci.\ U.S.A.}
{\bf 87} (1990),
1441--1445.

\itemitem{[6]}
{W.~Y.~C.~Chen and J.~D.~Louck},
{
The factorial Schur function},
{\sl J.\ Math.\ Phys.\ }
{\bf 34}
(1993),
4144--4160.

\itemitem {[7]} {P. H. Dondi and P. D. Jarvis}, {Diagram and
superfield techniques in the classical superalgebras},
{\sl J. Phys. A} {\bf 14} (1981), 547--563. 

\itemitem{[8]} {I. M. Gessel and G. X. Viennot}, {
Binomial determinants, paths, and hook length formulae},
{\sl Adv. Math.} {\bf 58} (1985), 300--321.

\itemitem{[9]} {I. Goulden and C. Greene,} A new tableau representation
for supersymmetric Schur functions, {\sl J. Algebra.} {\bf 170} (1994),
687--703.

\itemitem{[10]} {I. P. Goulden and A. M. Hamel,}
Shift operators and factorial symmetric functions, University of Waterloo 
Preprint CORR 92--10, April 1, 1992, to appear in 
{\it J. Comb. Theor. A.}

\itemitem{[11]}
{R. Howe},
{Remarks on classical invariant theory}, {\sl Trans. AMS}
{\bf 313}
(1989),
539--570.

\itemitem{[12]}
{R. Howe and T. Umeda},
{The Capelli identity, the double commutant theorem,
and multiplicity-free actions},
{\sl Math. Ann.}
{\bf 290}
(1991),
569--619.

\itemitem{[13]} {J. van der Jeugt, J. W. B. Hughes, R. C. King and
J. Thierry-Mieg}, {Character formulae for irreducible modules of the Lie
superalgebra $sl(m|n)$}, {\sl J. Math. Phys.} {\bf 31} (1990), 2278--2304.

\itemitem{[14]} {A.-A. A. Jucys}, {Symmetric polynomials and the center of
the symmetric group ring}, {\sl Reports Math. Phys.} {\bf 5} (1974),
107--112.

\itemitem{[15]} {V. G. Kac}, {Lie superalgebras}, {\sl Adv. Math.}
{\bf 26} (1977), 8--96.

\itemitem{[16]} {V. G. Kac}, {Representations of classical Lie
superalgebras}, {\sl in}
``Differential Geometry Methods in Mathematical Physics
II", (K. Bleuer, H. R. Petry, A. Reetz, Eds.),
Lecture Notes in Math., Vol. 676, pp. 597--626.
Springer-Verlag, Berlin/Heidelberg/New York, 1978.

\itemitem{[17]} {A. Lascoux}, {Letter to the author.}

\itemitem{[18]} {A. Lascoux}, {Classes de Chern des vari\'et\'es
de drapeaux}, {\sl Comptes Rendus Acad. Sci. Paris, S\'er. I}
{\bf 295} (1982), 393--398.

\itemitem{[19]} {D. E. Littlewood,} ``The theory of group characters and 
matrix representations of groups", 2nd edition, Clarendon Press,  Oxford,
1950.

\itemitem{[20]} {I. G. Macdonald,} 
``Symmetric functions and Hall polynomials", 2nd edition,
Oxford University Press, 1995.

\itemitem{[21]} {I. G. Macdonald},
{Schur functions: theme and variations}, {\sl in}
{``Actes 28-e S\'eminaire Lotharingien",\/} pp. 5--39.
{Publ. I.R.M.A. Strasbourg,\/} 1992, 498/S--27. 

\itemitem{[22]} {I. G. Macdonald},
``Notes on Schubert polynomials",
Publ. LACIM, Univ. du Qu\'ebec \`a Montr\'eal, 1991.

\itemitem{[23]} {N. Metropolis, G. Nicoletti, G. C. Rota},
A new class of symmetric functions, {\sl in} ``Mathematical Analysis
and Applications," pp. 563--575. Advances in Mathematics
Supplementary Studies, Vol. {\bf 7B}, Academic Press, New York, 1981.

\itemitem{[24]}
A. Molev, M. Nazarov,
{Capelli identities for classical groups},
Preprint University of Wales 95-21, Swansea, 1995.

\itemitem{[25]}
{A. I. Molev, M. L. Nazarov and G. I. Olshanski\u\i},
{Yangians and classical Lie algebras},
{\sl Uspekhi Matem. Nauk} {\bf 51}, (1996) no. 2, 27--104.

\itemitem{[26]} {G. E. Murphy}, {A new construction of Young's
seminormal representation of the symmetric group}, {\sl J. Algebra}
{\bf 69} (1981), 287--291.

\itemitem{[27]}
{M. L. Nazarov},
{Quantum Berezinian and the classical Capelli identity}, 
{\sl Lett. Math. Phys.}
{\bf 21}
(1991),
123--131.

\itemitem{[28]}
{M. L. Nazarov}, {Yangians and Capelli identities}, Preprint, 1995,
q-alg/9601027.

\itemitem{[29]}
A.~Okounkov,
{Quantum immanants and higher Capelli identities},
Preprint, 1995, q-alg/9602028.

\itemitem{[30]}
A.~Okounkov,
{Young basis, Wick formula and higher Capelli identities},
Preprint, 1995, q-alg/9602027.

\itemitem{[31]}
A.~Okounkov and G.~Olshanski\u\i,
{Shifted Schur functions},
Preprint, 1995, q-alg/9605042.

\itemitem{[32]}
{G. I. Olshanski\u\i},
{Quasi-symmetric functions and factorial Schur functions},
Preprint, 1995.

\itemitem {[33]}{P. Pragacz}, {Algebro-geometric applications of Schur
$S$- and $Q$-polynomials}, {\sl in} ``Topics in Invariant Theory,"
Seminaire d'Algebre Paul Dubriel et Marie-Paule Malliavin,
Proceedings, Lecture Notes in Math., Vol. 1478, pp. 130--191.
Springer-Verlag, New York/Berlin, 1991.

\itemitem {[34]}{P. Pragacz and A. Thorup}, {On a Jacobi--Trudi
identity for supersymmetric polynomials},
{\sl Adv. Math.} {\bf 95} (1992), 8--17.

\itemitem {[35]} {J. B. Remmel}, The combinatorics of $(k,l)$-hook
Schur functions, {\sl Contemp. Math.} {\bf 34} (1984), 253--287.

\itemitem {[36]} {B. Sagan}, ``The symmetric group: representations,
combinatorial algorithms, and symmetric functions", Wadsworth \& Brooks  
/ Cole mathematics series, Pacific Grove, CA, 1991.

\itemitem {[37]} {S. Sahi}, {The spectrum of certain invariant
differential operators associated to a Hermitian symmetric space},
{\sl in} ``Lie theory and geometry", Progr. Math. {\bf 123},
Boston: Birkhauser 1994, 569--576.

\itemitem {[38]} {M. Scheunert}, {Casimir elements of Lie superalgebras},
{\sl in} ``Differential Geometry Methods in Mathematical Physics",
pp. 115--124,
Reidel, Dordrecht, 1984.

\itemitem {[39]} {R. P. Stanley}, Unimodality and Lie superalgebras,
{\sl Stud. Appl. Math.} {\bf 72} (1985), 263--281.

\itemitem {[40]} {J. R. Stembridge}, {A characterization of supersymmetric
polynomials}, {\sl J. Algebra} {\bf 95} (1985), 439--444.

\enddocument